\documentclass[12pt]{amsart}
\usepackage{amsfonts,amssymb,amsmath,amscd,mathtools,amsthm,xcolor,mathrsfs,subfiles,tikz-cd,tikz}
\usepackage[backref,colorlinks, linkcolor=blue, citecolor=magenta]{hyperref}
\usepackage{latexsym}

\theoremstyle{plain}
\newtheorem{theorem}{Theorem}[section]
\newtheorem{proposition}[theorem]{Proposition}
\newtheorem{corollary}[theorem]{Corollary}
\newtheorem{lemma}[theorem]{Lemma}
\newtheorem{example}[theorem]{Example} 
\numberwithin{equation}{section}

\theoremstyle{definition}
\newtheorem{definition}[theorem]{Definition}

\theoremstyle{remark}
\newtheorem{remark}[theorem]{Remark}

\def\N{\mathbb{N}}
\def\Z{{\mathbb{Z}}}
\def\R{{\mathbb{R}}}

\def\Hom{{\rm Hom}}
\def\F{\mathcal{F}}
\def\G{\mathcal{G}}
\def\O{\mathcal{O}}
\def\P{\mathcal{P}}
\def\S{\mathcal{S}}
\def\m{\mathfrak{m}}
\def\p{\mathfrak{p}}
\DeclareMathOperator{\spec}{{\rm Spec}}
\newcommand{\smod}[1]{{\rm\bf SMod}\text{-}#1}
\newcommand{\SRing}[1]{{\rm\bf SRing}(#1)}
\newcommand{\SFgp}[1]{{\rm\bf SFgp}(#1)}
\newcommand{\omod}[1]{{\rm\bf SMod}\text{-}\O_#1}
\newcommand{\Slfb}[1]{{\rm\bf SLfb}\text{-}\O_#1}
\newcommand{\Shom}[3][]{\overline{\mathcal{H}om}_{#1}\left(#2,#3\right)}
\title{The Serre-Swan Theorem in supergeometry}
\author{Archana S. Morye, Abhay Soman, and V. Devichandrika}
\address{School of Mathematics and Statistics, University of Hyderabad, Hyderabad-500046, India}
\email{asmsm@uohyd.ac.in,abhaysoman@uohyd.ac.in,18mmpp03@uohyd.ac.in}
\thanks{The first and the third named authors are supported by `Institute of Eminence University of Hyderabad' (UoH-IoE-RC5-22-003). A. Soman acknowledges support from SERB MATRICS (MTR/2022/00093).}
\date{\today}
\begin{document}

\begin{abstract}
	We show the analogue of the Serre-Swan theorem in a context of supergeometry. This theorem gives an equivalence of the category of locally free supersheaves of bounded rank over locally ringed superspace with the category of finitely generated super projective modules over its coordinate superring, under the assumptions that every locally free supersheaf is generated by global sections and it is acyclic.
\end{abstract}
\maketitle

\section{Introduction}
A celebrated result of J.-P. Serre gives one-to-one correspondence between algebraic vector bundles over affine variety over algebraically closed field and finitely generated projective modules over its coordinate ring (see \cite{Serre}). Later in \cite{Swan}, R. G. Swan showed one-to-one correspondence between topological vector bundles over compact Hausdorff space $X$ and finitely generated projective modules over the ring of continuous real-valued functions on $X$. In this paper, Swan gave several nontrivial examples of projective modules using the correspondence. Results of Serre and Swan was later generalized to locally ringed spaces by Morye in \cite{Archana}.

The main result in this article is the following analogue of the Serre-Swan theorem for locally ringed superspaces.
\begin{theorem}\label{intro-main-theorem}
	Let $(X,\O_X)$ be a locally ringed superspace and $A=\Gamma(X,\O_X)$. Assume that every locally free supersheaf of bounded rank is acyclic and generated by finitely many global sections. Then the functor $\S$ defines an equivalence of categories from the category $\SFgp{A}$ of finitely generated superprojective right $A$-supermodules to the category $\Slfb{X}$ of locally free supersheaves of bounded rank over $X$.
\end{theorem}

The proof of the above theorem follows the proof given in \cite{Archana}, \cite{Morye-2013} and we use \cite{Westra} for the required supercommutative algebra. Perhaps the Serre-Swan theorem in supergeometry is known to experts in the area. However, we did not find a reference for this result.


This article is organized as follows. In Section \ref{prelim}, we recall basic definitions and some results related to supercommutative rings and supermodules. We define locally ringed superspaces and supersheaves in Section \ref{sec-basics-supersheaves}. Furthermore, we show in this section that the sheaf associated to a presheaf of superrings is again a sheaf of superrings, and the stalk at a point is a superring. In Section \ref{sec-category-sheaves-of-supermodules}, we define basic notions related to sheaves of supermodules, including sheaves generated by global sections and locally free supersheaves. The notions and constructions developed in this section are similar to the standard constructions in sheaf theory. We developed these notions in the context of supergeometry. In Section \ref{sec-adjoint-pair}, we start with a ringed superspace $(X,\O_{X})$, and we define a presheaf $\P(M)$ associated to $\O_X(X)$-supermodule $M$ by $\P(M)(U)=M\otimes_{\O_X(X)}\O_X(U)$. We consider $\S(M)$ to be the sheaf associated to $\P(M)$. This defines a functor $\S$ from the category of $\O_X(X)$-supermodules to the category of $\O_X$-supersheaves. We show that $\S$ and the global section functor $\Gamma$ are adjoint pairs and we study their properties.
In the penultimate section \ref{sec-super-serre-swan}, we prove Theorem \ref{intro-main-theorem}. In the last section we give examples of affine superschemes and smooth supermanifolds for which Theorem \ref{intro-main-theorem} holds.
\section{Preliminaries}\label{prelim}
\subsection{Definitions and Notations}We begin with some basic definitions and notions. For more details we refer to \cite{Bruzzo91}, \cite{Bruzzo} and \cite{Manin}. 

Let $\Z_2=\{0,1\}$, the field with two elements. A $\Z_2$-graded ring with unity is called a \emph{superring}. We denote the \emph{degree} of a homogeneous element $r$ in a superring by $|r|$. The degree $0$ elements of  a superring $A$ forms a ring and it is called the \emph{even component} of $A$, and the set of degree $1$ elements in $A$ is called the \emph{odd component} of $A$. A superring $A$ is said to be \emph{supercommutative} if it satisfies the condition $rs=(-1)^{|r||s|}sr$ for any homogeneous elements $r,s\in A$. We denote by $J_A$ the ideal generated by odd elements of a superring $A$. Note that, $J_A$ is a homogeneous ideal and that every element of $J_A$ is nilpotent. The quotient ring $A/J_A$, which is again a superring, is called the \emph{Bosonic reduction}.  A supercommutative ring $A$ is said to be a \emph{local superring} if it has a unique homogeneous maximal ideal. A \emph{homomorphism of superrings} $\varphi\colon A\to B$ is a ring homomorphism such that $\varphi(A_\delta)\subseteq B_{\delta}$ (where $\delta=0,1$).

We give a typical example of superring.
\begin{example}\label{typical-example}
	Consider a polynomial ring over noncommuting variables $\beta_1,\ldots,\beta_n$:
	\[\R[\beta_1,\beta_2,\ldots,\beta_n]=\left\{\sum a_{k_1k_2\cdots k_n}\beta_1^{k_1}\cdots\beta_n^{k_n}:a_{k_1k_2\cdots k_n}\in\R\text{ and } k_i\in\Z_{\geq 0}\right\}\]
	satisfying the following conditions.
	\[\beta_i^2=0\quad\text{and}\quad\beta_i\beta_j=-\beta_j\beta_i\quad\text{for all}\quad 1\leq i,j\leq n.\]
	This is a supercommutative ring. 
	We denote $\R[\beta_1,\beta_2,\ldots,\beta_n]$ by $\Lambda$. The even (resp. odd) component of $\Lambda$ is denoted by $\Lambda_0$ (resp., $\Lambda_1$). We denote by $J_\Lambda$ the ideal generated by $\beta_1,\beta_2,\ldots,\beta_n$. Let  $\pi\colon\Lambda\to\Lambda\big/J_\Lambda$ be the canonical homomorphism of superrings. We remark that $\Lambda\big/J_\Lambda$ and $\R$ are isomorphic superrings with $\R_0=\R$ and $\R_1=0$. By abuse of notation, we consider the canonical homomorphism of superrings $\pi\colon\Lambda\to\R$.
\end{example} 

Let $A$ be a supercommutative ring. A right $\Z_2$-graded $A$-module (resp., left $\Z_2$-graded $A$-module) $M$ is called a \emph{supermodule over $A$} or an \emph{$A$-supermodule}, i.e., $M=M_0\bigoplus M_1$ and  $M_{\delta}A_\epsilon\subseteq M_{\epsilon+\delta}$ (resp., $A_\epsilon M_\delta\subseteq M_{\epsilon+\delta}$) for $\epsilon,\delta\in\Z_2$. We denote by $\Pi M$ a right $A$-supermodule with $\left(\Pi M\right)_0=M_1$ and $\left(\Pi M\right)_1=M_0$. The right action of $A$ on $\Pi M$ is same as the right action of $A$ on $M$. Let $M,N$ be right $A$-supermodules. A \emph{homomorphism of supermodules} $u\colon M\to N$ is an $A$-linear map such that $u\left(M_\delta\right)\subseteq N_\delta$ for $\delta\in\Z_2$. The set of all homomorphisms of $A$-supermodules from $M$ to $N$ is denoted by $\Hom_A(M,N)$. The set of all $A$-module homomorphisms (not necessarily $\Z_2$-graded) from $M$ to $N$ is denoted by $\overline{\Hom}_A\left(M,N\right)$. We remark that $\overline{\Hom}_A\left(M,N\right)$ is an $A$-supermodule, and the even component of $\overline{\Hom}_A\left(M,N\right)$ is precisely $\Hom_A(M,N)$. For a supercommutative ring $A$, the collection of right $A$-supermodules with $A$-supermodule homomorphisms forms a category which is denoted by $\smod{A}$.

Suppose that $N$ is a right $A$-supermodule. We give a left $A$-supermodule structure on $N$ as follows.
\[a\cdot n\coloneqq (-1)^{|a||n|}\,na\]
Let $M,N$ be right $A$-supermodules. We consider $N$ as a left $A$-supermodule as described above. The tensor product $M\otimes_AN$ is defined in an usual way.
Note that for $a\in A$, $m\in M$, and $n\in N$ we have the following relations.
\[(m\otimes n)a=ma\otimes n=m\otimes an=m\otimes(-1)^{|a||n|}na=(-1)^{|a||n|}\,m\otimes na\]

A right \emph{free $A$-supermodule} of rank $p|q$ is an $A$-supermodule isomorphic to $A^p\bigoplus\left(\Pi A\right)^q$, where $A$ is considered as a right $A$-supermodule over itself. The module $A^p\bigoplus\left(\Pi A\right)^q$ is denoted by $A^{p|q}$. We say that a right $A$-supermodule $M$ is \emph{finitely presented} if we have the following exact sequence of $A$-supermodules for some $p,q,r,s\in\N$.
\[A^{r|s}\to A^{p|q}\to M\to 0\]

A right $A$-supermodule $P$	is said to be \emph{projective supermodule} if $P$ is a direct summand of a right free $A$-supermodule. We refer to \cite[Section 6.2]{Westra} and \cite[Section 3]{MSV}.        

\subsection{Super-analogues of a few results from commutative algebra}
In this subsection, we state and prove super-analogues of some results from commutative algebra for completeness. We follow \cite{Westra}. Throughout this subsection we assume that $A$ is a supercommutative ring and all morphisms are in the category $\smod{A}$.

\begin{lemma}\label{prime-ideal-implies-graded-prime}
	Let $A$ be a supercommutative ring and let $\mathfrak{p}$ be a prime ideal (not necessarily assumed to be homogeneous). Then, $\mathfrak{p}$ is homogeneous and contains $J_A$.
\end{lemma}
\begin{proof}
	As odd elements are nilpotent, $A_1\subseteq\mathfrak{p}$ and hence, $J_A\subseteq\mathfrak{p}$.
	For any $a\in\mathfrak{p}$ there exists a unique $a_i\in A_i$ for $i=0,1$ such that $a=a_0+a_1$. As $a_1\in J_A\subseteq\mathfrak{p}$, we also get $a_0\in\mathfrak{p}$. This shows that $\mathfrak{p}$ is a homogeneous ideal.
\end{proof}

\begin{lemma}
	Let $A=A_0\oplus A_1$ be a supercommutative ring. We have the following one to one correspondences.
	\begin{align}
		\left\{\text{prime ideals of }A_0\right\}&\longleftrightarrow\left\{\text{homogeneous prime ideals of }A\right\}\label{even-prime-correspondence}\\
		\left\{\text{maximal ideals of }A_0\right\}&\longleftrightarrow\left\{\text{homogeneous maximal ideals of }A\right\}\label{even-maximal-correspondence}
	\end{align}
\end{lemma}
\begin{proof}
	We prove the correspondence for prime ideals. The correspondence for maximal ideals may be proved in the same way. By Lemma \ref{prime-ideal-implies-graded-prime}, a prime ideal $\p$ of $A$ is of the form $\p_0\oplus A_1$, where $\p_0$ is an ideal of $A_0$. We define the map $A_0\big/\p_0\to A\big/\p$ given by $a_0+\p_0\mapsto a_0+\p$. This is a homomorphism of superrings. It is easy to see that this map is an isomorphism. Hence the result.
\end{proof}

\begin{lemma}[Free supermodules and tensor product]\label{free-supermodules-tensor-product}
	Let $M$ be a right $A$-supermodule and $A^p\bigoplus\left(\Pi A\right)^q$ be a free right $A$-supermodule of rank $p|q$. Then we have the following isomorphism of right $A$-supermodules.
	\[\left(A^p\bigoplus\left(\Pi A\right)^q\right)\otimes_AM\simeq M^p\bigoplus\left(\Pi M\right)^q\]
\end{lemma}

\begin{proof}
	Define a map $\varphi\colon \left(A^p\bigoplus\left(\Pi A\right)^q\right)\otimes_AM\longrightarrow M^p\bigoplus\left(\Pi M\right)^q$ as follows.
	\[\left((x_1,\ldots,x_p,y_1,\ldots,y_q),m\right)\mapsto\left(m\cdot x_1,\ldots,m\cdot x_p,m\cdot y_1,\ldots,m\cdot y_q\right)\]
	The map $\varphi$ is an $A$-supermodule homomorphism.
	
	Assume that for $1\leq i\leq p+q$, $e_i=(\delta_{ij})_{j=1}^{p+q}\in A^p\bigoplus\left(\Pi A\right)^q$. We define an $A$-supermodule homomorphism 
	\begin{align*}
		\psi\colon M^p\bigoplus\left(\Pi M\right)^q&\longrightarrow\left(A^p\bigoplus\left(\Pi A\right)^q\right)\otimes_AM\text{ by }\\
		(m_1,m_2,\ldots,m_{p+q})&\mapsto\sum_{i=1} ^{p+q}e_i\otimes m_i
	\end{align*}

	Furthermore, maps $\varphi$ and $\psi$ are inverses of each other.
\end{proof}

The following lemma shows that the inverse limit commutes with tensoring with a finitely generated superprojective modules.

\begin{lemma}\label{inverse-limit-commutes-with-finitely-presented}
	Let $I$ be a nonempty poset. Let $A$ be a supercommutative ring and let $M_i$ ($i\in I$) be a right $A$-supermodules. Suppose that $P$ is a finitely generated superprojective right $A$-supermodule. Then we get the following isomorphism of right $A$-supermodules.
	\[P\otimes\varprojlim M_i\simeq\varprojlim (P\otimes M_i)\]
\end{lemma}
\begin{proof}
	Consider the canonical morphism $\gamma\colon P\otimes_A\varprojlim M_i\to\varprojlim(P\otimes_AM_i)$ given by $x\otimes(m_i)\mapsto (x\otimes m_i)$. We first show that $\gamma$ is injective.
	
	As $P$ is a finitely generated superprojective module, it is finitely presented (see, \cite[Proposition 6.5.7]{Westra}). A super-analogue of the result from \cite[Chapter I, Exercise 9]{Bourbaki-commutative} shows that the canonical morphism $\widetilde{\gamma}\colon P\otimes_A\prod M_i\to\prod (P\otimes_AM_i)$ given by $x\otimes(m_i)\to (x\otimes m_i)$ is bijective. We consider the following diagram.
	\[\begin{tikzcd}
		0 & {P\otimes_A\varprojlim M_i} & {P\otimes_A\prod M_i} \\
		0 & {\varprojlim(P\otimes_AM_i)} & {\prod(P\otimes_AM_i)}
		\arrow[from=1-1, to=1-2]
		\arrow[hook, from=1-2, to=1-3]
		\arrow["\gamma"', from=1-2, to=2-2]
		\arrow["{\widetilde{\gamma}}", from=1-3, to=2-3]
		\arrow[from=2-1, to=2-2]
		\arrow[hook, from=2-2, to=2-3]
	\end{tikzcd}\]
	By \cite[Proposition 6.2.5]{Westra}, $P$ is a flat $A$-supermodule, and hence, the first horizontal row in the above diagram is exact. Furthermore, as $\gamma$ is the restriction of $\widetilde{\gamma}$, the above diagram is commutative. Therefore, $\gamma$ is injective.
	
	We now show that $\gamma$ is surjective. As $P$ is finitely generated superprojective, it is a direct summand of a free $A$-supermodule $F$ of finite rank, say $F\simeq P\oplus K$, (see \cite[Lemma 6.2.2]{Westra}). Therefore, for each $i$, $F\otimes_AM_i\simeq (P\otimes_AM_i)\oplus(K\otimes_AM_i)$. 
	By \cite[Chapter IX, Section 2, Theorem 1]{MacLane}, finite colimits commute with limits. Thus, we get the following canonical isomorphisms.
	\[F\otimes_A\varprojlim M_i\stackrel{\alpha}{\simeq} \varprojlim(F\otimes_AM_i)\simeq\varprojlim(P\otimes_AM_i)\oplus\varprojlim(K\otimes_AM_i)\] In particular, \[0\to\varprojlim(K\otimes_AM_i)\to \varprojlim(F\otimes_AM_i)\xrightarrow{}\varprojlim(P\otimes_AM_i)\to 0\] is an exact sequence. Consider the following commutative diagram with exact rows.
	\[\begin{tikzcd}[sep=small]
		0 & {K\otimes_A\varprojlim M_i} & {F\otimes_A\varprojlim M_i} & {P\otimes_A\varprojlim M_i} & 0 \\
		0 & {\varprojlim(K\otimes_AM_i)} & {\varprojlim(F\otimes_AM_i)} & {\varprojlim(P\otimes_AM_i)} & 0
		\arrow[from=1-1, to=1-2]
		\arrow[from=1-2, to=1-3]
		\arrow[from=1-2, to=2-2]
		\arrow["{}", from=1-3, to=1-4]
		\arrow["\alpha"', from=1-3, to=2-3]
		\arrow[from=1-4, to=1-5]
		\arrow["\gamma", from=1-4, to=2-4]
		\arrow[from=2-1, to=2-2]
		\arrow[from=2-2, to=2-3]
		\arrow[from=2-3, to=2-4]
		\arrow[from=2-4, to=2-5]
	\end{tikzcd}\]
	As $\alpha$ is an isomorphism, $\gamma$ is surjective. Therefore, $\gamma$ is a required isomorphism.
\end{proof}

We recall from \cite[Theorem 6.4.7]{Westra} the super-analogue of Cayley-Hamilton theorem.

\begin{theorem}[Cayley-Hamilton]\label{Cayley-Hamilton}
	Let $A$ be a supercommutative ring and $I$ a homogeneous ideal of $A$. Let $M$ be a finitely generated right $A$-supermodule, and $\varphi\colon M\to M$ a supermodule homomorphism such that $\varphi(M)\subseteq MI$. Then there exists a monic polynomial $p(x)=x^n+a_1x^{n-1}+\cdots+a_n\in A_0[x]$ such that $p(\varphi)=0$ with $a_i\in I^i$.
\end{theorem}

\begin{corollary}\cite[Corollary 6.4.8]{Westra}\label{corollary-to-Cayley-Hamilton}
	Let $A$ be a supercommutative ring and let $I$ be a homogeneous ideal of $A$. Suppose that $M$ is a finitely generated right $A$-supermodule such that $M=MI$. Then there exists an even element $a\in A_0$ such that $M(1+a)=0$.
\end{corollary}

\begin{corollary}\label{surjective-module-hom-is-isomorphism}
	Let $A$ be a supercommutative ring and let $M$ be a finitely generated right $A$-supermodule. If $\varphi\colon M\to M$ is a surjective supermodule homomorphism then, $\varphi$ is an isomorphism.
\end{corollary}
\begin{proof}
	Consider the polynomial ring in one variable over $A$, $A[x]$. The ring $A[x]$ is a supercommutative ring. Indeed, the even (resp. odd) part of $A[x]$ can be considered to be polynomials with coefficients from $A_0$ (resp. $A_1$). We define a right $A[x]$-module structure on $M$ by considering $m\cdot x\coloneqq\varphi(m)$. As $\varphi$ is surjective, for $I=(x)$ we have $M=MI$. By Corollary \ref{corollary-to-Cayley-Hamilton}, there exists $p(x)\in A_0[x]$ such that $M(1+p(x))=0$, i.e., for any $m\in M$ we get $m+p(\varphi)(m)=0$. Hence, $\varphi$ is injective and therefore, $\varphi$ is an isomorphism.
\end{proof}
\section{Basics of supersheaves}\label{sec-basics-supersheaves}
In this section, we collect some basic results related to supersheaves. For general notions related to locally ringed spaces and sheaves we refer to relevant chapters from \cite{Hartshorne}. Throughout this section all superrings are assumed to be supercommutative. We begin with the following definition.
\begin{definition}\label{def-presheaf-of-superring}
	Let $X$ be a topological space. A presheaf $\mathcal{F}$ on $X$ is said to be a \emph{presheaf of superrings} if for every open subset $U$ of $X$, $\mathcal{F}(U)$ is a superring and whenever $V\subseteq U$, the restriction map $\rho_V^U\colon\mathcal{F}(U)\to\mathcal{F}(V)$ is a homomorphism of superrings.
	A \emph{sheaf of superrings} on a topological space $X$ is a presheaf of superrings which is a sheaf.
\end{definition}

\begin{lemma}\label{lemma-germ-is-superring}
	Let $X$ be a topological space and let $\F$ be a presheaf of superrings on $X$. For every $x\in X$, the stalk of $\F$ at $x$, $\F_x$ is a superring.
\end{lemma}
\begin{proof}
	As $\F$ is a presheaf of superrings, $\F_x$ is a ring for every $x\in X$. We show that 
	\[\F_x=\F_{x,0}\bigoplus\F_{x,1},\]
	where,
	\[\F_{x,0}=\left\{\left(s_0\right)_x:\text{where }s_x\in\F_x\quad\text{and}\quad s=s_0+s_1\in\F(U)\text{ for some open }U\right\}\]
	\[\F_{x,1}=\left\{\left(s_1\right)_x:\text{where }s_x\in\F_x\quad\text{and}\quad s=s_0+s_1\in\F(U)\text{ for some open }U\right\}.\]
	
	Let $s_x\in\F_x$, where $s\in\F(U)$ for some $U$. Suppose that $s=s_0+s_1\in\F(U)$ and $t=t_0+t_1\in\F(V)$, with $s_0\in\F(U)$, $t_0\in\F(V)$ are even and $s_1\in\F(U)$, $t_1\in\F(V)$ are odd elements, are such that $s_x=t_x\in\F_x$. We claim that $\left(s_0\right)_x=\left(t_0\right)_x$ and $\left(s_1\right)_x=\left(t_1\right)_x$. Suppose that $W$ is an open subset of $U\cap V$ containing $x$ such that $s|_W=t|_W$. As the restriction maps are homomorphisms of superrings, we get the following.
	\[s|_W=s_0|_W+s_1|_W\in \F(W)\quad\text{and}\quad t|_W=t_0|_W+t_1|_W\in\F(W)\]
	As $\F(W)$ is a superring and $s_0|_W,t_0|_W\in\F(W)$ (resp., $s_1|_W,t_1|_W\in\F(W)$) are even (resp., odd elements) components of $s|_W$ and $t|_W$, we get that $s_0|_W=t_0|_W$ and $s_1|_W=t_1|_W$. Therefore,
	\begin{equation}\label{eq-germ-is-superring}
		\left(s_0\right)_x=\left(t_0\right)_x\in\F_x\quad\text{and}\quad \left(s_1\right)_x=\left(t_1\right)_x\in\F_x.
	\end{equation}
	
	We thus proved that
	$\F_x=\F_{x,0}\bigoplus\F_{x,1}$.
	
	Let $(s_1)_x,(t_1)_x\in\F_{x,1}$. Consider an open subset $W$ of $X$ containing $x$ such that $s_1|_Wt_1|_W\in\F(W)$. As $\F(W)$ is a superring and restriction is a homomorphism of superrings, both $s_1|_W$ and $t_1|_W$ are odd elements of $\F(W)$ and hence, $s_1|_Wt_1|_W=-t_1|_Ws_1|_W\in\F(W)$. Thus, $(s_1)_x(t_1)_x=-(t_1)_x(s_1)_x\in\F_x$. Similarly, other homogeneous cases will follow. This shows that $\F_x$ is a superring.
\end{proof}

\begin{lemma}[Gluing of homogeneous sections]\label{Gluing-of-homogeneous-sections}
	Let $\F$ be a sheaf of superrings on a topological space $X$. Suppose that $U=\bigcup_{i\in I}U_i$ is an open covering of an open subset $U\subseteq X$. Let $s_i\in\F(U_i)$ be sections of odd degrees (resp., even degrees) such that for every $i,j\in I$ the following condition holds.
	\[s_i|_{U_i\cap U_j}=s_j|_{U_i\cap U_j}\]
	Then, there exists a unique section $s\in\F(U)$ of \emph{odd degree} (resp., \emph{even degree}) such that $s|_{U_i}=s_i$ ($i\in I$).  
\end{lemma}
\begin{proof}
	There exists a unique section
	$s\in\F(U)$ such that $s|_{U_i}=s_i$ for every $i\in I$, because $\F$ is assumed to be a sheaf. We need to show that the degree of $s\in\F(U)$ is odd (resp., even). Suppose that $s=s_0+s_1\in\F(U)$. As a restriction is a homomorphism of superrings we have the following. \[s|_{U_i}=s_0|_{U_i}+s_1|_{U_i}=s_i\in\F(U_i)\quad\text{for each}\quad i\in I\] By the assumption, each $s_i\in\F(U_i)$ is of odd degree (resp., even degree) and $\F(U_i)$ is a superring. Therefore, $s_0|_{U_i}=0$ (resp., $s_1|_{U_i}=0$) for every $i\in I$. As $\F$ is a sheaf, $s_0=0\in\F(U)$ (resp., $s_1=0\in\F(U)$), in particular, $s=s_0+s_1=s_1\in\F(U)$  ($s=s_0+s_1=s_0\in\F(U)$), i.e., $s\in\F(U)$ is of odd degree (resp., even degree).
\end{proof}

\begin{corollary}\label{locally-homoge-globally-homege}
	Let $\F$ be a sheaf of superrings on a topological space $X$. Let $U\subseteq X$ be an open subset. If $s\in\F(U)$ is such that $s_x\in\F_x$ is odd (resp., even) for all $x\in U$, then $s\in\F(U)$ is odd (resp., even).
\end{corollary}
\begin{proof}
	As $\F(U)$ is  a superring, $s=s_0+s_1$. By \eqref{eq-germ-is-superring}, $s_x=(s_0)_x+(s_1)_x\in\F_x$, for every $x\in U$. Hence by hypothesis, $s_x=(s_1)_x$ (resp., $s_x=(s_0)_x$). As $\F$ is a sheaf we get that $s=s_1\in\F(U)_1$ (resp., $s=s_0\in\F(U)_0$).
\end{proof}

\begin{lemma}\label{sheaf-associated-is-supersheaf}
	Let $X$ be a topological space and $\F$ be a presheaf of supperrings on $X$. The sheaf $\widetilde{\F}$ associated with $\F$ is a sheaf of superrings.
\end{lemma}
\begin{proof}
	For an open subset $U\subseteq X$, we consider  $\widetilde{\F}(U)\subseteq \prod_{x\in U}\F_x$ to be the collection of all $\theta=\left(\theta(x)_{x\in U}\right)$ such that for every $y\in U$ there exists an open neighborhood  $V$ of $y$ such that $s_x=\theta(x)$ for all $x\in V$.
	We show that $\widetilde{\F}(U)$ is a superring. We consider the following sets
	\[\widetilde{\F}(U)_0=\left\{(\theta(x))\in \widetilde{\F}(U):\theta(x)\in\F_{x,0}\right\}\]
	and
	\[\widetilde{\F}(U)_1=\left\{(\theta(x))\in \widetilde{\F}(U):\theta(x)\in\F_{x,1}\right\}.\]
	We claim that $\widetilde{\F}(U)=\widetilde{\F}(U)_0\bigoplus\widetilde{\F}(U)_1$. Let $\theta=(\theta(x))\in\widetilde{\F}(U)$. For every $x\in U$, $\theta(x)=\theta(x)_0+\theta(x)_1\in\F_x=\F_{x,0}\bigoplus\F_{x,1}$ (see Lemma \ref{lemma-germ-is-superring}). Put $\theta_0=(\theta(x)_0)$ and $\theta_1=(\theta(x)_1)$. We show that $\theta_0,\theta_1\in\widetilde{\F}(U)$. For any $y\in U$, there exists an open neighborhood $V\subseteq U$ of $y$ and $s=s_0+s_1\in\F(V)$ such that $\theta(x)=s_x$. By the proof of Lemma \ref{lemma-germ-is-superring}, $s_x=(s_0)_x+(s_1)_x$ and hence $\theta(x)_i=(s_i)_x$ (for $i=1,2$). This shows that $\theta_i\in\widetilde{\F}(U)$ for $i=1,2$. As in the proof of Lemma \ref{lemma-germ-is-superring}, $\widetilde{\F}(U)$ satisfies supercommutative condition. Hence the result.
\end{proof}

\begin{definition}[Category of sheaves of superrings]\label{def-category-sheaves-of-superring}
	The category of sheaves of superrings $\SRing{X}$ on a topological space on $X$ is defined as follows.
	\begin{enumerate}
		\item Objects of $\SRing{X}$ are sheaves of superrings on $X$.
		\item If $\F$ and $\G$ are sheaves of superrings then,  a morphism of sheaves $\varphi\colon\F\to\G$ is called a \emph{morphism of sheaves of superrings} if for every open subset $U\subseteq X$, $\varphi(U)\colon\F(U)\to\G(U)$ is a homomorphism of superrings. 
	\end{enumerate}
\end{definition}

\begin{example}\label{supersheaf-of-continous-functions}
	We consider superring $\Lambda$ as given in Example \ref{typical-example}. Note that $\Lambda$ is a $2^n$-dimensional vector space over $\R$. A topology on $\Lambda$ is induced from $\R^{2^n}$. For an open subset $U\subseteq X$, let $\mathcal{C}_\Lambda(U)=\left\{f\colon U\to \Lambda:f\text{ is continuous}\right\}$. It is easy to see that $\mathcal{C}_\Lambda$ is a sheaf of rings with usual restriction maps. In fact, $\mathcal{C}_\Lambda(U)$ is a superring with the following even and odd components.
	\begin{align*}
		\mathcal{C}_\Lambda(U)_0&=\left\{f\in\mathcal{C}_\Lambda(U):f(U)\subseteq \Lambda_0\right\}\\
		\mathcal{C}_\Lambda(U)_1&=\left\{f\in\mathcal{C}_\Lambda(U):f(U)\subseteq \Lambda_1\right\}
	\end{align*}
	Note that for any open subsets $V\subseteq U$, the restriction map $\mathcal{C}_\Lambda(U)\to\mathcal{C}_\Lambda(V)$ is a homogeneous map. Therefore, $\mathcal{C}_\Lambda$ is a sheaf of superrings.
\end{example}

Let $\F$ be a presheaf and $\G$ a sheaf. For a given family of morphisms $\F_x\to\G_x$, for each $x$, the following proposition gives a criterion to extend this family to a unique morphism $\F\to\G$.
This is one of the key result we require to prove the main theorem.
\begin{proposition}\label{local-morphism-to-global-morphism}
	Let $\F$ be a presheaf and $\G$ a sheaf on a topological space $X$. For every $x\in X$, suppose that there exist morphisms $\alpha^x\colon\F_x\to\G_x$. Further assume that for every open subset $U\subseteq X$ and a section $s\in\F(U)$ there exists, for all $x\in U$, an open neighborhood $V_x\subseteq U$ of $x$ and $t\in\G(V_x)$ such that $t_y=\alpha^y(s_y)$ for every $y\in V_x$. Then there exists a unique morphism of presheaves $\alpha\colon\F\to\G$ such that $\alpha_x=\alpha^x$, where $\alpha_x\colon\F_x\to\G_x$ is the induced stalk level map from $\alpha$.
\end{proposition}
\begin{proof}
	Define $\alpha_U\colon\F(U)\to\G(U)$ as follows. Let $s\in\F(U)$. By the assumption, there exists an open covering $\{V_i\}$ of $U$ and $t_i\in\G(V_i)$ such that $(t_i)_y=\alpha^y(s_y)$. For every $y\in V_i\cap V_j$, $(t_i)_y=\alpha^y(s_y)=(t_j)_y$. Therefore, $t_i|_{V_i\cap V_j}=t_j|_{V_i\cap V_j}$. As $\G$ is a sheaf, there exists $t\in\G(U)$ such that $t_y=(t_i)_y$ for $y\in V_i$. We define $\alpha_U\colon\F(U)\to\G(U)$ by $s\mapsto t$. Note that this map is well-defined. Indeed, if $\{W_j\}$ is an open cover of $U$ and $t_j'\in\G(W_j)$ such that $(t'_j)_y=\alpha^y(s_y)$. We obtain $t'\in\G(U)$ such that $t'_y=(t_j')_y$ for $y\in W_j$. Note that for every $x\in U$, $t_x=\alpha^x(s_x)=t_x'$. Since $\G$ is a sheaf, $t=t'$.
	
	By similar arguments one can show that $\alpha$ is a morphism of presheaves.
	
\end{proof}

We remark that if we assume in Proposition \ref{local-morphism-to-global-morphism} that $\F$ is a presheaf of superrings (resp., supermodules) and $\G$ is a sheaf of superrings (resp., supermodules), then by arguments similar to Corollary \ref{locally-homoge-globally-homege}, $\alpha$ is a morphism of presheaves of superrings (resp., supermodules).

\section{Category of sheaves of supermodules}\label{sec-category-sheaves-of-supermodules}
In this section, we gather some definitions and notations related to sheaves of supermodules that will be required in the subsequent sections. All the definitions are standard and can be found in \cite{Hartshorne}, however we give definitions in context of superrings. We begin with the following definitions.

\begin{definition}[Ringed superspace]\label{def-ringed-superspace}
	A pair $(X,\mathcal{O}_X)$ is said to be a \emph{ringed superspace} if $\mathcal{O}_X$ is a sheaf of superrings on $X$.
\end{definition}

\begin{definition}[Locally ringed superspace]\label{locally-ringed-superspace}
	A ringed superspace $(X,\mathcal{O}_X)$ is said to be \emph{locally ringed superspace} if for a given $x\in X$, $\mathcal{O}_{X,x}$ is a local superring.
\end{definition}

\begin{example}
	A ringed space $(X,\mathcal{C}_\Lambda)$ defined in Example \ref{supersheaf-of-continous-functions} is an example of a locally ringed superspace. For any $x\in X$, the maximal ideal of $\mathcal{C}_{\Lambda,x}$ is
	\[\m_x=\left\{s_x\in\mathcal{C}_{\Lambda,x}:\left(\pi\circ s\right)(x)=0\right\}\] where, $\pi\colon\Lambda\to\R$ is the canonical map as defined in Example \ref{typical-example}
	
\end{example}

\begin{example}[Affine superscheme]\label{Affine-superscheme}
	Let $A$ be a superring and let $J_A$ be an ideal generated by odd elements of $A$. By \eqref{even-prime-correspondence}, there is a one-one correspondence between prime ideals of $A$ and prime ideals of $A_0$, i.e.,  ${\rm Spec}(A)\leftrightarrow{\rm Spec}(A_0)\leftrightarrow{\rm Spec}(A/J_A)$. Put $X={\rm Spec}(A_0)$. The affine superscheme $(X,\O_X)$ is defined in such a way that for $f\in A_0$, $\O_X\left(D(f)\right)=A_f$. This is an example of locally ringed superspace. For more details refer to  \cite[Section 5.4]{Westra}.
\end{example}

We now recall a few definitions related to sheaves of supermodules.

\begin{definition}[Sheaf of supermodule]\label{def-Ox-supermodule}
	Let $(X,\mathcal{O}_X)$ be a ringed superspace. A sheaf $\F$ on $X$ is called a \emph{sheaf of $\mathcal{O}_X$-supemodule} if for every open subset $U\subseteq X$, $\F(U)$ is a right $\mathcal{O}_X(U)$-supermodule, and the following conditions are satisfied for every open subsets $V\subseteq U$.
	\begin{enumerate}
		\item Restriction map $\F(U)\to\F(V)$ is homomorphisms of superrings.
		\item The following diagram commutes for $V\subseteq U$.
		\[\begin{tikzcd}
			\F(U)\times\O_X(U) && \F(U) \\
			\F(V)\times\O_X(V) && \F(V)
			\arrow[from=1-1, to=1-3]
			\arrow[from=1-1, to=2-1]
			\arrow[from=2-1, to=2-3]
			\arrow[from=1-3, to=2-3]
		\end{tikzcd}\]
		Here, vertical arrows represent restriction maps and horizontal arrows represent the scalar multiplication.
	\end{enumerate} 
\end{definition}

\begin{definition}[Category of right $\O_X$-supermodules]\label{def-category-of-right-Ox-supermodules}
	Let $(X,\O_X)$ be a ringed superspace. The category of right $\O_X$-supermodules, $\omod{X}$ is defined as follows.
	\begin{enumerate}
		\item Objects are sheaves of right $\O_X$-supermodules.
		\item Let $\F,\G$ be right $\O_X$-supermodules. A homomorphism $u\colon\F\to\G$ is a homomorphism of sheaves such that for any
		open subset $U\subseteq X$, $u_U\colon\F(U)\to\G(U)$ is $\O_X(U)$-supermodule homomorphism.  
	\end{enumerate}
\end{definition}

\begin{remark}\label{Ox-supermodules-is-abelian-category}
	The category of $\O_X$-supermodules is an abelian category.
\end{remark}

We denote by $\overline{\Hom}_{\O_X}(\F,\G)$ by the set of all $\O_X$-module homomorphisms from $\F$ to $\G$. Note that $\overline{\Hom}_{\O_X}(\F,\G)$ is a $\Gamma(X,\O_X)$-supermodule.  Furthermore, we denote by $\Hom_{\O_X}(\F,\G)$ the set of homomorphism of $\O_X$-supermodules, which is the even component of $\overline{\Hom}_{\O_X}(\F,\G)$.

\begin{definition}[Sheaf of homomorphisms of $\mathcal{O}_X$-supermodules]\label{Hom-supersheaf}
	Let $\F$ and $\G$ be right $\mathcal{O}_X$-supermodules. The sheaf of homomorphism of $\O_X$-supermodules is denoted by $\mathcal{H}om_{\O_X}(\F,\G)$. For any open subset $U\subseteq X$, we define $\mathcal{H}om_{\O_X}(\F,\G)(U)=\Hom_{\O_X|_U}(\F|_U,\G|_U)$. The restriction maps are the usual restriction maps. 
\end{definition}

\begin{definition}[Supersheaf of homomorphisms of $\mathcal{O}_X$-supermodules]\label{def-supersheaf-of-complete-hom}
	Let $\left(X,\mathcal{O}_X\right)$ be a ringed superspace. Let $\F$ and $\G$ be right $\O_X$-supermodules. A sheaf of all $\O_X$-module homomorphisms from $\F$ to $\G$ is denoted by $\Shom[\O_X]{\F}{\G}$.
\end{definition}

We remark that the even component of $\Shom[\O_X]{\F}{\G}$ is $\mathcal{H}om_{\O_X}(\F,\G)$.

We now state a few definitions and prove some results related to $\O_X$-supermodules that will be used in the proof of the main theorem \ref{super-serre-swan}.

\begin{definition}[Global section functor]\label{def-Gamma-functor}
	Let $(X,\O_X)$ be a locally ringed space. We recall the left exact \emph{global section functor} \[\Gamma(X,\cdot)\colon\omod{X}\to\smod{A}\] given by $\F\mapsto\F(X)$.
\end{definition}

\begin{definition}[Parity swapping functor]\label{def-parity-swapping-functor-for-O-mod}
	Let $(X,\O_X)$ be a ringed superspace.  The \emph{parity swapping functor} $\Pi\colon\omod{X}\to\omod{X}$ is defined as follows. For $\F\in{\rm Ob}\left(\omod{X}\right)$, we associate $\Pi\F\in{\rm Ob}\left(\omod{X}\right)$ such that $(\Pi\F)(U)=\Pi\left(\F(U)\right)$, for any open subset $U\subseteq X$ and $\Pi\left(\F(U)\right)$ is defined as in Section \ref{prelim} for a supermodule $\F(U)$. It is easy to verify that $\Pi\F$ is a supersheaf.
\end{definition}

\begin{example}
	Let $(X,\O_X)$ be a ringed superspace and let $I,J$ be nonempty indexing sets. We define $\O_X$-module $\O_X^{(I|J)}$ as follows. For every open subset $U\subseteq X$,
	\[\O_X^{(I|J)}(U)\coloneqq\left(\bigoplus_I\O_X(U)\right)\bigoplus\left(\bigoplus_J\Pi\O_X(U)\right).\]
	If $|I|=r<\infty$ and $|J|=s<\infty$ then, we denote $\O_X^{(I|J)}$ by $\O_X^{r|s}$. 
\end{example}

\begin{definition}[Locally free supersheaf]\label{def-locally-free-sheaf}
	Let $(X,\O_X)$ be a ringed superspace. A right $\O_X$-module $\F$ is said to be \emph{locally free supersheaf} if for any $x\in X$ there exists an open neighborhood $U\subseteq X$ of $x$ and nonempty sets $I,J$ (may depend on $x$) such that \[\F|_U\simeq\O_X^{(I|J)}|_U\quad\text{as}\quad\O_X|_U\text{-supermodule}.\]
	
	A locally free right $\O_X$-supermodule $\F$ is said to be \emph{locally free of finite rank} if  $|I|<\infty$ and $|J|<\infty$.
\end{definition}

\begin{definition}[Quasi-coherent supersheaf]\label{def-Quasi-coherent-supersheaf}
	Let $(X,\O_X)$ be a ringed superspace. A  right $\O_X$-supermodule $\F$ is said to be \emph{Quasi-coherent} $\O_X$-supermodule if for every $x\in X$ there exists an open neighborhood $U\subseteq X$ such that for some sets $I_1,I_2,J_1,J_2$ (may depend on $x$) we have the following exact sequence.
	\[(\O_X|_U)^{(I_1|J_1)}\to(\O_X|_U)^{(I_2|J_2)}\to\F|_U\to 0\]
\end{definition}

\begin{definition}[Supersheaf of finite type]\label{def-finite-type}
	Let $(X,\O_X)$ be a ringed superspace. A right $\O_X$-supermodule $\F$ is said to be of \emph{finite type} if for every $x\in X$, there exists an open neighborhood $U$ of $x$ such that $\F|_U$ is generated by a finite family of sections. 
\end{definition}

\begin{remark}\label{quotient of finite type}
	Let $(X,\O_X)$ be a locally ringed superspace, and let $\F,\G$ be right $\O_X$-supersheaves. If $\F$ is of finite type and if $\varphi\colon\F\to\G$ is a surjective morphism, then $\G$ is of finite type. In particular, every quotient supersheaf of finite type is of finite type.
\end{remark}

\begin{definition}[Finitely presented supersheaf]\label{def-finitely-presented-supersheaf}
	Let $(X,\O_X)$ be a ringed superspace. A  right $\O_X$-supermodule $\F$ is said to be \emph{finitely presented} $\O_X$-supermodule if for every $x\in X$ there exists an open neighborhood $U\subseteq X$ such that for some $r,s,p,q\in\N$ (may depend on $x$) we have the following exact sequence.
	\[(\O_X|_U)^{r|s}\to(\O_X|_U)^{p|q}\to\F|_U\to 0\]
\end{definition}

\begin{definition}[Sheaf generated by global sections]\label{sheaf-generated-by-global-sections}
	Let $(X,\O_X)$ be ringed superspace and let $\F$ be a right $\O_X$-supermodule. We say that $\F$ is \emph{generated by global sections} if there is a family of global sections $s_i\in\Gamma(X,\F)$, where $i\in I$ for some indexing set $I$, such that for every $x\in X$, $\F_x$ is generated by $(s_i)_x$ ($i\in I$) as an $\O_{X,x}$-supermodule.
\end{definition}

\begin{remark}\label{sheaf-generated-by-homogeneous-global-sections}
	We keep notations of Definition \ref{sheaf-generated-by-global-sections}. Suppose that $s_i=s_{i0}+s_{i1}$, where $s_{i0}$ is an even element and $s_{i1}$ is an odd element. Then, $\F$ is generated by homogeneous global sections  $\{s_{i0},s_{i1}\}_{i\in I}$. \emph{Throughout this article we assume that if a sheaf is generated by global sections then these sections will be assumed to be homogeneous}.
\end{remark}

\begin{proposition}\label{support-is-closed-for-finite-type}
	Let $(X,\O_X)$ be a locally ringed superspace. If $\F$ is of finite type, then the support of $\F$, ${\rm Supp}(\F)=\{x\in X:\F_x\neq 0\}$ is a closed subset of $X$.
\end{proposition}
\begin{proof}
	Let $x\in X\setminus {\rm Supp}(\F)$, i.e., $\F_x=0$. As $\F$ is of finite type, there exists a neighborhood $U$ of $x$ and  $t_1,t_2,\ldots,t_n\in\F(U)$ such that $\{(t_1)_y,\ldots,(t_n)_y\}$ generates $\F_y$ for every $y\in U$. In particular, $(t_i)_x=0$ for each $i$. Thus, there exists $V\subseteq U$ such that $t_i|_V=0$ for all $i$. Therefore, $\F_y=0$ for all $y\in V$, i.e., $V\subseteq X\setminus{\rm Supp}(\F)$. Hence the result.
\end{proof}

\begin{corollary}\label{surjective at a point implies in neighborhood for finite type}
	Let $(X,\O_X)$ be a locally ringed superspace, and let $\F,\G$ be right $\O_X$-supersheaves of finite type. Suppose that $\varphi\colon\F\to\G$ is a morphism of supersheaves. If $\varphi_x\colon\F_x\to\G_x$ is surjective at $x\in X$ then, there exists an open neighborhood $U$ of $x$ such that $\varphi|_U\colon\F|_U\to\G|_U$ is surjective.
\end{corollary}
\begin{proof}
	As $\G$ is of finite type, the supersheaf $\G\big/{\rm Im}(\varphi)$ is of finite type. By Proposition \ref{support-is-closed-for-finite-type}, ${\rm Supp}\left(\G\big/{\rm Im}(\varphi)\right)$ is a closed subset of $X$. Furthermore, as $\varphi_x$ is surjective, $x\in X\setminus{\rm Supp}\left(\G\big/{\rm Im}(\varphi)\right)$. Hence, there exists an open neighborhood $U$ of $x$ such that $\varphi|_U\colon\F|_U\to\G|_U$ is surjective.
\end{proof}

\begin{proposition}\label{isomorphism-of-Hom-supersheaf}
	Let $(X,\O_X)$ be a locally ringed superspace. If $\F$ is a finitely presented $\O_X$-supersheaf and $\G$ an $\O_X$-supersheaf then for every $x\in X$ there exists $\O_{X,x}$-supermodule isomorphism of superrings
	\[\phi_x\colon\left(\overline{\mathcal{H}om}_{\O_X}\left(\F,\G\right)\right)_x\longrightarrow\overline{\Hom}_{\O_{X,x}}\left(\F_x,\G_x\right).\]
	In particular, $\left({\mathcal{H}om}_{\O_X}\left(\F,\G\right)\right)_{x}\simeq {\Hom}_{\O_{X,x}}\left(\F_x,\G_x\right)$ as abelian groups.
\end{proposition}
\begin{proof}
	We have $\left(\left(\overline{\mathcal{H}om}_{\O_X}\left(\F,\G\right)\right)_x\right)_0=\left({\mathcal{H}om}_{\O_X}\left(\F,\G\right)\right)_{x}$ by Corollary \ref{locally-homoge-globally-homege}. Note that the map $\phi_x$ is defined as follows. For $\alpha\in\left(\overline{\mathcal{H}om}_{\O_X}\left(\F,\G\right)\right)_x$ consider an open neighborhood $U$ of $x$ and a morphism $u\colon\F|_U\to\G|_U$ such that $\alpha=u_x$. Then we define $\phi_x(\alpha)\coloneqq u_x$. By \cite[Proposition 1.13, page no. 213]{Archana}, $\phi_x$ is an isomorphism of $\O_{X,x}$-modules. As $\phi_x$ is a morphism of supermodules, the even component of $\left(\overline{\mathcal{H}om}_{\O_X}\left(\F,\G\right)\right)_x$ is isomorphic to the even component of $\overline{\Hom}_{\O_{X,x}}\left(\F_x,\G_x\right)$. Hence the result.
\end{proof}

\begin{corollary}\label{finitely-presented-germs-to-local}
	Let $\F$ and $\G$ be an $\O_X$-supermodules of finite presentation, and let $x\in X$. Let $f\colon\F_x\to\G_x$ be an isomorphism of $\O_{X,x}$-supermodules. Then, there exists an open neighborhood $U$ of $x$, and an isomorphism of $\O_X|_U$-supermodules $u\colon\F|_U\to\G|_U$ such that $u_x=f$. 
\end{corollary}
\begin{proof}
	Let $g\colon\G_x\to\F_x$ be such that $f\circ g=1_{\G_x}$ and $g\circ f=1_{\F_x}$. By Proposition \ref{isomorphism-of-Hom-supersheaf}, there exists open neighborhood $U\subseteq X$ of $x$ and   $u\in\mathcal{H}om_{\O_X|_U}(\F|_U,\G|_U)$ such that $u_x=f$. Similarly, there exists $v\colon\G|_U\to\F|_U$ such that $v_x=g$. As $u_x\circ v_x=1_{\G_x}$ and $v_x\circ u_x=1_{\F_x}$, there exists a neighborhood $V\subseteq U$ of $x$ such that $u|_V\circ v|_V=1_{\G|_V}$, and $v|_V\circ u|_V=1_{\F|_V}$. Hence the result.
\end{proof}

We get the following result as a consequence of Corollary \ref{finitely-presented-germs-to-local}.

\begin{corollary}\label{finitely-presented-free-germs-to-locally-free}
	Let $(X,\O_X)$ be a ringed superspace and let $\F$ be an $\O_X$-supermodule of finite presentation. Let $x\in X$ and let $\F_x$ be a free $\O_{X,x}$-supermodule. Then, there exists an open neighborhood $U$ of $x$ and natural numbers $r$ and $s$ such that $\F|_U\simeq\O_X^{r|s}|_U$.	In particular, if $\F_x$ is a free $\O_{X,x}$-supermodule for every $x\in X$, then $\F$ is locally free.
\end{corollary}

\section{Properties of an adjoint pair $\Gamma$ and $\S$}\label{sec-adjoint-pair}

The super-Serre-Swan theorem gives an equivalence of categories of locally free supersheves of bounded rank with finitely generated  superprojective modules. This equivalence is given by the global section functor and the functor $\S$ that we define below.

\begin{definition}[$\S$ functor]\label{def-tensor-functor}
	Let $(X,\O_X)$ be a ringed superspace and $A=\Gamma(X,\O_X)$. Let $M$ be a right $A$-supermodule. For an open subset $U\subseteq X$, the restriction map $\O_X(X)\to\O_X(U)$ makes $\O_X(U)$ a right $A$-supermodule, viz., the scalar multiplication is defined by $s\cdot a=s(a|_U)$.
	
	We define a presheaf $\mathcal{P}(M)$ associated to $M$ of $\O_X$-supermodules as follows. For an open subset $U\subseteq X$, 
	\begin{equation}\label{def-eq-tensor-functor}
		\P(M)(U)=M\otimes_A\O_X(U).
	\end{equation}
	
	We denote the sheaf associated with $\P(M)$ by $\S(M)$. Note that $\S(M)$ is a sheaf of right $\O_X$-supermodules.
	We define the following functor.
	\begin{align*}
		\S\colon\smod{A}&\to\omod{X}\\
		M&\longmapsto\S(M)\\
		u\colon M\to N&\longmapsto u\otimes 1_{\O_X}
	\end{align*}
\end{definition}

\begin{lemma}\label{finitely-presented-implies-PM-is-sheaf}
	We keep the notations of Definition \ref{def-tensor-functor}. If $M$ is a finitely generated superprojective right $A$-supermodule, then $\P(M)$ is a sheaf.
\end{lemma}
\begin{proof}
	By Lemma \ref{inverse-limit-commutes-with-finitely-presented}, \[M\otimes_A\varprojlim_{V\subseteq U}\O_X(V)\simeq\varprojlim_{V\subseteq U}\left(M\otimes_A\O_X(V)\right)=\S(M)(U).\]
	As $\O_X$ is a sheaf, we have $\O_X(U)=\displaystyle\varprojlim_{V\subseteq U}\O_X(V)$. Hence, $\P(M)(U)=\S(M)(U)$ and thus, $\P(M)$ is a sheaf.
\end{proof}

\begin{lemma}\label{S(P) is locally free of bounded rank}
	Let $(X,\O_X)$ be locally ringed superspace and let $A=\Gamma(X,\O_X)$. If $P$ is a finitely generated  superprojective $A$-supermodule then, $\S(P)$ is locally free $\O_X$-module of finite rank.
\end{lemma}
\begin{proof}
	By \cite[Lemma 6.5.7]{Westra}, $P$ is finitely presented. Suppose that there is the following exact sequence of right $A$-supermodules.
	\[A^{r|s}\to A^{p|q}\to P\to 0\]
	As $\S$ is a right exact functor, we get the following exact sequence of $\O_X$-supermodules.
	\[\S\left(A^{r|s}\right)\to\S\left(A^{p|q}\right)\to\S(P)\to 0\]
	We have $\S\left(A^{r|s}\right)\simeq\O_X^{r|s}$ and $\S\left(A^{p|q}\right)\simeq\O_X^{p|q}$ (see Lemma \ref{free-supermodules-tensor-product}) and hence the following sequence is exact.
	\[\O_X^{r|s}\to\O_X^{p|q}\to\S(P)\to 0\]
	For any $x\in X$, $\S(P)_x\simeq P\otimes_A\O_{X,x}$, which is a superprojective $\O_{X,x}$-supermodule. By \cite[Lemma 6.4.6]{Westra}, $\S(P)_x$ is a free $\O_{X,x}$-supermodule of finite rank. By Corollary \ref{finitely-presented-free-germs-to-locally-free}, $\S(P)$ is locally free of finite rank.
\end{proof}

\begin{proposition}\label{left-right-exactness-Gamma-tensor}
	Let $(X,\O_X)$ be a locally ringed superspace. The global section functor $\Gamma(X,\cdot)$ and the functor $\S$ is an adjoint pair.
\end{proposition}
\begin{proof}
	Suppose that $A=\Gamma(X,\O_X)$, $M$ a right $A$-supemodule, and $\F$ a right $\O_X$-supermodule. We define a map 
	\begin{equation}\label{def-lambda-map}
		\lambda_{\F,M}\colon\Hom_A(M,\Gamma(X,\F))\to\Hom_{\O_X}(\S(M),\F).
	\end{equation}
	For $u\in\Hom_A(M,\Gamma(X,\F))$, we define $\lambda_u\colon\P(M)\to\F$ as follows. For an open subset $U\subseteq X$, and $m\in M$ and $s\in\O_X(U)$ we consider 
	\[\lambda_u(m\otimes s)=u(m)|_U\cdot s.\]
	By the universal property of sheaf associated with a presheaf, $\lambda_u$ will induce an $\O_X$-supermodule homomorphism $\lambda_{\F,M}(u)\colon\S(M)\to\F$. 
	Let $g\colon M'\to M$ be an $A$-supermodule homomorphism, and let $\psi\colon\F\to\F'$ be an $\O_X$-supermodule homomorphism. We consider, $\Hom(g,\Gamma(X,\psi))(u)$ to be the composition of the following $A$-supermodule homomorphisms.
	\[\begin{tikzcd}
		M' & M & \Gamma(X,\F) & \Gamma(X,\F')
		\arrow["g", from=1-1, to=1-2]
		\arrow["u", from=1-2, to=1-3]
		\arrow["\psi", from=1-3, to=1-4]
		\arrow["{\Hom(g,\Gamma(X,\psi))(u)}"', bend right=25, dashed, from=1-1, to=1-4]
	\end{tikzcd}\]
	Similarly, $\Hom(\S(g),\psi)(\alpha)$ is the composition of the following sheaf maps.
	\[\begin{tikzcd}
		\S(M') & \S(M) & \F & \F'
		\arrow["\S(g)", from=1-1, to=1-2]
		\arrow["\alpha", from=1-2, to=1-3]
		\arrow["\psi", from=1-3, to=1-4]
		\arrow["{\Hom(\S(g),\psi)(\alpha)}"', bend right=25, dashed, from=1-1, to=1-4]
	\end{tikzcd}\]
	To show that $\Gamma(X,\cdot)$ and $\S$ is an adjoint pair we prove that the following diagram commutes.
	\[\begin{tikzcd}
		\Hom_A(M,\Gamma(X,\F)) && \Hom_{\O_X}(\S(M),\F) \\
		\Hom_A(M',\Gamma(X,\F')) && \Hom_{\O_X}(\S(M'),\F')
		\arrow["\lambda_{\F,M}",from=1-1, to=1-3]
		\arrow["{\Hom(g,\Gamma(X,\psi))}"',from=1-1, to=2-1]
		\arrow["\lambda_{\F',M'}"',from=2-1, to=2-3]
		\arrow["{\Hom(\S(g),\psi)}",from=1-3, to=2-3]
	\end{tikzcd}\]
	Let $U$ be an open subset of $X$, $m'\in M'$, and $s\in\O_X(U)$. We get the following.
	\begin{align*}
		\left(\Hom(\S(g),\psi)\circ\lambda_{\F,M}(u)\right)_U(m'\otimes s)&=\left(\psi\circ\lambda_{\F,M}(u)\circ\S(g)\right)_U(m'\otimes s)\\
		&=\left(\psi\circ\lambda_{\F,M}(u)\circ\S(g)\right)_X(m')|_U\cdot s\\
		&=\psi(u(g(m')))|_U\cdot s
	\end{align*}
	On the other hand,
	\begin{align*}
		\left(\lambda_{\F',M'}\circ\Hom(g,\Gamma(X,\F))(u)\right)_U(m'\otimes s)&=\left(\lambda_{\F',M'}(\psi\circ u\circ g)\right)_U(m'\otimes s)\\
		&=\psi\left(u(g(m'))\right)|_U\cdot s
	\end{align*}
	This shows the commutativity of the diagram, and hence $\Gamma(X,\cdot)$ and $\S$ is an adjoint pair.
\end{proof}

\begin{proposition}\label{Gamma-fully-faithful}
	Let $(X,\O_X)$ be ringed superspace and let $\mathcal{C}$ be a full abelian subcategory of $\omod{X}$ such that $\O_X\in{\rm Ob}(\mathcal{C})$. Suppose that every sheaf in $\mathcal{C}$ is generated by global sections. Then $\Gamma(X,\cdot)\colon\mathcal{C}\to\omod{X}$ is fully faithful.
\end{proposition}
\begin{proof}
	Let $A=\Gamma(X,\O_X)$ and $\F,\G$ in $\mathcal{C}$. We need to show the following map
	\[\phi_{\F,\G}\colon\Hom_{\O_X}(\F,\G)\to\Hom_A\left(\Gamma(X,\F),\Gamma(X,\G)\right)\quad\text{given by}\quad u\mapsto u_X\]
	is bijective. Let $\F$ be generated by homogeneous global sections $\{s_i\}_{i\in I}\cup\{t_j\}_{j\in J}$, where $s_i$ are even and $t_j$ are odd (see Remark \ref{sheaf-generated-by-homogeneous-global-sections}).
	
	We first prove the injectivity of $\phi_{\F,\G}$. Suppose that $u\colon\F\to\G$ is an $\O_X$-supermodule morphism such that $u_X=0$. Hence $u_X(s_i)=0$ and $u_X(t_j)=0$, and therefore, $u_x\colon\F_x\to\G_x$ is the zero map for any $x\in X$. Thus, $u$ is a zero morphism of supersheaves.
	
	We now show that $\phi_{\F,\G}$ is surjective. Let $\alpha\colon\Gamma(X,\F)\to\Gamma(X,\G)$ be a right $A$-supermodule homomorphism. For every $x\in X$, we define the map $\alpha^x\colon\F_x\to\G_x$ as follows.
	\[\alpha^x\left(\sum (s_i)_x\cdot a_i+\sum (t_j)_x\cdot b_j\right)=\sum (\alpha(s_i))_x\cdot a_i+\sum (\alpha(t_j))_x\cdot b_j\]
	where, $a_i,b_j\in\O_{X,x}$ and $a_i=0$ (resp. $b_j=0$) for all but finitely many $i$ (resp., $j$). To show that $\alpha^x$ is a well-defined homomorphism it is enough to show that if $\sum_{i\in I_1} (s_i)_x\cdot a_i+\sum_{j\in J_1} (t_j)_x\cdot b_j=0$ then, $\sum_{i\in I_1} (\alpha(s_i))_x\cdot a_i+\sum_{j\in J_1} (\alpha(t_j))_x\cdot b_j=0$, where $I_1,J_1\subseteq I$ and $|I_1|=n<\infty, |J_1|=m<\infty$. In that vein, we consider a homomorphism of $\O_X$-supermodules $\varphi\colon\O_X^{(I_1|J_1)}\to\F$. For any open subset $U\subseteq X$, the map $\varphi_U\colon\O_X^{(I_1|J_1)}(U)\to\F(U)$ is defined in the following way.
	\[\varphi_U\left((c_i)_{j\in I_1},(d_j)_{j\in J_1}\right)=\sum_{i\in I_1}s_i|_U\cdot c_i+\sum_{j\in J_1}t_j|_U\cdot d_j\] 
	Let $\mathcal{K}=\ker\varphi$. Note that $\left((a_i)_{i\in I_1},(b_j)_{j\in J_1}\right)\in\mathcal{K}_x$. Since $\mathcal{K}\in{\rm Ob}(\mathcal{C})$ there exists global sections $(g_k,h_k)$ ($1\leq k\leq r$) of $\mathcal{K}$, where $g_{k}=(g_{ki})_{i\in I_1}\in A^{I_1}$ and $h_k=(h_{kj})_{j\in J_1}\in(\Pi A)^{J_1}$, such that \[\left((a_i)_{i\in I_1},(b_j)_{j\in J_1}\right)=\sum_{k=1}^r(g_k,h_k)_x\cdot e_k,\quad\text{for some } e_k\in\O_{X,x}.\]
	
	As $(g_k,h_k)\in\Gamma(X,\mathcal{K})$, $\sum_{i\in I_1}s_i\cdot g_{ki}+\sum_{j\in J_1}t_j\cdot h_{kj}=0$ for all $1\leq k\leq r$. We have the following.
	\begin{align*}
		&\sum_{i\in I_1} (\alpha(s_i))_x\cdot a_i+\sum_{j\in J_1} (\alpha(t_j))_x\cdot b_j\\
		&=\sum_{i\in I_1} (\alpha(s_i))_x\cdot \left(\sum_{k=1}^r (g_{ki})_x\cdot e_k\right)+\sum_{j\in J_1} (\alpha(t_j))_x\cdot \left(\sum_{k=1}^r (h_{kj})_x\cdot e_k\right)\\
		&=\sum_{k=1}^{r}\left(\sum_{i\in I_1} (\alpha(s_i))_x\cdot (g_{ki})_x\right)\cdot e_k+\sum_{k=1}^{r}\left(\sum_{j\in J_1} (\alpha(t_j))_x\cdot (h_{kj})_x\right)\cdot e_k\\
		&=\sum_{k=1}^r\left(\sum_{i\in I_1}\left(\alpha(s_i)\cdot g_{ki}\right)_x\right)\cdot e_k+\sum_{k=1}^r\left(\sum_{j\in J_1}\left(\alpha(t_j)\cdot h_{kj}\right)_x\right)\cdot e_k
	\end{align*}
	As $\alpha$ is a homomorphism of supermodules, we get that
	\begin{align*}
		&\sum_{k=1}^r\sum_{i\in I_1}\left(\alpha(s_i\cdot g_{ki})\right)_x\cdot e_k+\sum_{i=1}^k\sum_{j\in J_1}\left(\alpha(t_j\cdot h_{kj})\right)_x\cdot e_k\\
		&=\sum_{k=1}^r\left(\alpha\left(\sum_{i\in I_1}s_i\cdot g_{ki}+\sum_{j\in J_1}t_j\cdot h_{kj}\right)\right)_x\cdot e_k=0				
	\end{align*}
	This shows that $\alpha^x$ is well-defined.
	Furthermore, as $\alpha$ is a homomorphism of $A$-supermodules, $\alpha^x$ is  homomorphism of $\O_{X,x}$-supermodules. Indeed,
	\begin{align*}
		&\alpha^x\left(\sum_{i\in I_1} (s_i)_x\cdot a_i+\sum_{j\in J_1}(t_j)_x\cdot b_j\right)\\
		&=\alpha^x\left(\sum(s_i)_x\cdot (a_{i0}+a_{i1})+\sum(t_j)_x\cdot (b_{j0}+b_{j1})\right)\\
		&=\sum\alpha(s_i)_x\cdot (a_{i0}+a_{i1})+\sum\alpha(t_j)_x\cdot (b_{j0}+b_{j1})\\
		&=\sum_{i,j}\left(\alpha(s_i)_x\cdot a_{i0}+\alpha(t_j)\cdot b_{j1}\right)+\sum_{i,j}\left(\alpha(s_i)_x\cdot a_{i0}+\alpha(t_j)\cdot b_{j0}\right).
	\end{align*}
	Note that the the even (resp., odd) component of $\sum_{i\in I_1} (s_i)_x\cdot a_i+\sum_{j\in J_1}t_j\cdot b_j$ is 
	\[\sum_{i,j}(s_i)_x\cdot a_{i0}+(t_j)_x\cdot b_{j1}\quad\left(\text{resp., }\sum_{i,j}(s_i)_x\cdot a_{i1}+(t_j)_x\cdot b_{j0}\right).\]
	
	The family $\{\alpha^x\}$ will give a unique morphism $u\colon\F\to\G$ of supersheaves as follows. 
	
	Let $U$ be an open subset of $X$ and $z\in\F(U)$. For $x\in U$, there exists even $s_{i_1},\ldots,s_{i_n}$ and odd $t_{j_1},\ldots,t_{j_m}$ such that $z_x=\sum_{k=1}^n(s_{i_k})_x\cdot a_k+\sum_{\ell=1}^m(t_{j_\ell})_x\cdot b_\ell$, where $a_k,b_\ell\in\O_{X,x}$. There exists an open subset $V\subseteq U$ containing $x$ and $f_k,g_\ell\in\O_X(V)$ such that $a_k=(f_k)_x$, $b_\ell=(g_\ell)_x$, and $z|_V=\sum s_{i_k}|_V\cdot f_k+\sum t_{\ell_k}|_V\cdot g_\ell$. Let $w=\sum \alpha(s_{i_k})|_V\cdot f_k+\sum\alpha(t_{\ell_k})|_V\cdot g_\ell\in\G(V)$. Then $\alpha^y(z_y)=w_y$ for all $y\in V$.    	
	We now apply Proposition \ref{local-morphism-to-global-morphism} to get a unique morphism $u\colon\F\to\G$ such that $\alpha^x=u_x$. Furthermore, $u_X=\alpha$. This shows the surjectivity of $\phi_{\F,\G}$.
\end{proof}

\begin{proposition}\label{key-prop-to-main-theorem}
	Let $(X,\O_X)$ be ringed superspace. Let $\mathcal{C}$ be a full abelian subcategory of $\omod{X}$ such that $\O_X\in{\rm Ob}(\mathcal{C})$, and every sheaf in $\mathcal{C}$ is generated by global sections. Then, every supersheaf $\mathcal{F}$ in $\mathcal{C}$ is isomorphic to $\mathcal{S}(\Gamma(X,\F))$.
\end{proposition}
\begin{proof}
	
	Let $M=\Gamma(X,\F)$.  We define a morphism of presheaves of right $\O_X$-supermodules $u':\mathcal{P}(M)\to \F$ as follows. For $U\subset X$ open, $(u')_U:M\otimes \mathcal{O}_X(U)\to \F(U)$ is given by $(u')_U(m\otimes f)=m|_U\cdot f$ for $m\in M$, and  $f\in\mathcal{O}_X(U)$.  By the universal property of sheafification,  we get a morphism of  supersheaves $u:\mathcal{S}(M)\to \F$ associated to $u'$.  We show that $u$ is an isomorphism by proving for all $x\in X$, $u_x:M\otimes \mathcal{O}_{X,x}\to\F_x$ is an isomorphism.    
	
	We first prove that $u_x$ is surjective. Since $\F$ is generated by global sections, for  a given $w\in\F_x$ there exists  even global section $s_i$ of $\F$ where, $1\leq i\leq n$, and odd global sections $t_j$ of $\F$ where, $1\leq j\leq m$, and $a_i,b_j\in\mathcal{O}_{X,x}$ such that $w=\sum_{i=1}^n (s_i)_x\cdot a_i+\sum_{j=1}^m (t_j)_x\cdot b_j$.   Then,  $u_x(\sum_{i=1}^n s_i\otimes a_i+\sum_{j=1}^m t_j\otimes b_j)=w$, hence $u_x$ is surjective. 
	
	Now we show that $u_x$ is injective. If  \[u_x\left(\sum_{i=1}^n s_i\otimes a_i+\sum_{j=1}^m t_j\otimes b_j\right)=\sum_{i=1}^n (s_i)_x\cdot a_i+\sum_{j=1}^m (t_j)_x\cdot b_j=0\] where $s_i$ are even and $t_j$ are odd global sections of $\F$, and $a_i,b_j\in \mathcal{O}_{X,x}$, then we prove $\sum_{i=1}^n s_i\otimes a_i+\sum_{i=1}^m t_j\otimes b_j=0$.   Consider a morphism of supersheaves $\phi:\mathcal{O}_X^{n|m}\to \F$, given by the family $\{s_i\}_{i=1}^n\cup\{t_j\}_{j=1}^m$. That is, for a given open subset $U$ of $X$, $\phi_U:\mathcal{O}_X^{n|m}(U)\to \F(U)$ is  defined by   $(\alpha_1,\ldots,\alpha_n,\beta_1,\ldots,\beta_m)\mapsto \sum_{i=1}^ns_i|_U\cdot \alpha_i+\sum_{j=1}^{m}t_j|_U\cdot\beta_j$, for $\alpha_i,\beta_j\in \mathcal{O}_{X}(U)$.  Note that $(a_1,\ldots,a_n,b_1,\ldots,b_m)\in \mathcal{K}_x$, where $\mathcal{K}=\ker(\phi)$.  Since $\mathcal{K}$ is in the category $\mathcal{C}$, $\mathcal{K}$ is generated by global sections.  Therefore, there exists $(g_k,h_k)=(g_{k1},\ldots,g_{kn},h_{k1},\ldots,h_{km})\in\Gamma(X,\mathcal{K})$ such that \[(a_1,\ldots,a_n,b_1,\ldots,b_m)=\sum_{k=1}^p (g_k)_x\cdot z_k+(h_k)_x\cdot w_k,\] for some $z_k,w_k\in\mathcal{O}_{X,x}$.  For $i=1,\ldots, n$, $a_i=\sum_{k=1}^p (g_{ki})_x\cdot z_k+(h_{ki})_x\cdot w_k$, $b_j=\sum_{k=1}^p (g_{kj})_x\cdot z_k+(h_{kj})_x\cdot w_k$ and $\sum_{i=1}^n s_i \cdot g_{\ell i}+\sum_{j=1}^{m}t_j\cdot h_{\ell j}=0$ for all $\ell=1,\ldots, p$. Thus,
	
	\begin{align*}
		&\sum_{i=1}^n s_i\otimes a_i+\sum_{j=1}^{m}t_j\otimes b_j\\
		& = \sum_{i=1}^n s_i \otimes \left(\sum_{k=1}^p (g_{ki})_x\cdot z_k+(h_{ki})_x\cdot w_k\right)\\
		&\qquad+\sum_{j=1}^{m}t_j\otimes\left(\sum_{k=1}^{p}(g_{kj})_x\cdot z_k+(h_{kl})_x\cdot w_k\right)\\
		&= \sum_{k=1}^p \left(\sum_{i=1}^n s_i\cdot g_{ki}+\sum_{j=1}^{m}t_j\cdot h_{kj}\right)\otimes z_k\\&\qquad+\sum_{k=1}^p \left(\sum_{i=1}^n s_i\cdot g_{ki}+\sum_{j=1}^{m}t_j\cdot h_{kj}\right)\otimes w_k\\
		&=0.
	\end{align*}
	Hence, $u_x$ is injective.  
\end{proof}

\section{Super Serre-Swan theorem}\label{sec-super-serre-swan}
In this section we prove the super-analogue of the Serre-Swan theorem. We begin with the following lemma.

\begin{lemma}\label{locally constant then locally free of finite rank}
	Let $(X, \mathcal{O}_X)$ be a locally ringed superspace and let $\F$ be a supersheaf of right $\O_X$-supermodule of finite type.  Then, $\F$ is locally free of finite rank if and only if $\F_x$ is free $\mathcal{O}_{X,x}$-supermodule and the rank of $\F$ is locally constant function.
\end{lemma}
\begin{proof}
	We only show that if $\F_x$ is a free $\O_{X,x}$-supermodule and the rank of $\F$ is locally constant function then, $\F$ is locally free of finite rank.  Let $x\in X$ and let $m|n$ be the rank of $\F_x$. By the hypothesis, there exists an open neighborhood $U$ of $x$ such that the rank of $\F_z$ is $m|n$ for any $z\in U$.  Thus, there is an isomorphism of $\O_{X,z}$-supermodules, $\psi^z\colon\F_z\to\O_{X,z}^{m|n}$ for every $z\in U$. We may assume that there exists even sections $s_1,s_2,\ldots,s_m\in\Gamma(U,\F)$ and odd sections $t_1,t_2,\ldots,t_n\in\Gamma(U,\F)$ such that $\mathfrak{B}=\{(s_1)_x,\ldots,(s_m)_x,(t_1)_x,\ldots,(t_n)_x\}$ is a basis of $\F_x$. Let $\varphi\colon\O_{X}|_U^{m|n}\to\F|_U$ be a homomorphism of $\O_X|_U$-supermodule defined $\mathfrak{B}$. Then, $\varphi_x\colon\O_{X,x}^{m|n}\to\F_x$ is an isomorphism of $\O_{X,x}$-supermodules. By Corollary \ref{surjective at a point implies in neighborhood for finite type}, there exists a neighborhood $V\subseteq U$ such that $\varphi|_V\colon\O_X|_V^{m|n}\to\F|_V$ is surjective morphism. For every $z\in V$, $\varphi_z\colon\O_{X,z}^{m|n}\to\F_z$ is surjective morphism. Hence, for each $z\in V$, $\psi_z\circ\varphi_z\colon\O_{X,z}^{m|n}\to\O_{X,z}^{m|n}$ is surjective. By Corollary \ref{corollary-to-Cayley-Hamilton}, $\psi_z\circ\varphi_z$ is an isomorphism for every $z\in V$. As $\psi_z$ is an isomorphism, so is $\varphi_z$. Hence, $\F$ is locally free of finite rank.
\end{proof}

\begin{proposition}\label{kernel-is-locally-free-of-locally-free-sheaves}
	Let $(X, \mathcal{O}_X)$ be a locally ringed superspace and $\F$, $\G$ locally free right $\mathcal{O}_X$-supermodules of finite rank.  If $u:\F\to\G$ is a sujective morphism, then $\ker(u)$ is also locally free of finite rank. 
\end{proposition}
\begin{proof}
	We get the following exact sequence of right $\O_X$-supermodules
	\[  0 \to \ker(u)\to \F\stackrel{u}{\to} \G\to 0. \]
	Since $\G$ is locally free the above sequence is locally split. Hence, there exist an open neighborhood $U$ of $x$, such that  $\ker(u)|_U\oplus \G|_U\simeq \F|_U$. This implies that $\ker(u)$ is an $\O_X$-supermodule of finite type.  Since $\F$ and $\G$ are locally free of finite rank, the rank of  $\ker(u)$ is locally constant.
	
	Since $\F_x, \G_x$ are free $\mathcal{O}_{X,x}$-supermodules, $\ker(u)_x$ is superprojective over local superring $\mathcal{O}_{X,x}$.   Hence, by \cite[Lemma 6.4.6]{Westra} $\ker(u)_x$ is free.  By the above Lemma \ref{locally constant then locally free of finite rank}, $\ker(u)$ is locally free of finite rank.    		   
\end{proof}

We denote the full subcategory of finitely generated superprojective $A$-supermodules by $\SFgp{A}$, and the subcategory of locally free $\O_X$-supermodules of bounded rank is denoted by $\Slfb{X}$.

\begin{proposition}\label{S is fully faithful}
	Let $(X,\O_X)$ be a locally ringed superspace and $A=\Gamma(X,\O_X)$. Then the functor $\S\colon{\bf Fgsp}(A)\to\Slfb{X}$ is fully faithful.
\end{proposition}

\begin{proof}
	By Lemma \ref{S(P) is locally free of bounded rank}, for every finitely generated superprojective right $A$-supermodule $P$, $\S(P)$ belongs to $\Slfb{X}$. 
	
	Suppose that $P$ and $Q$ are finitely generated superprojective $A$-supermodules. In particular, $P$ and $Q$ are finitely presented $A$-supermodules. We show that $\S\colon\Hom_A(P,Q)\to\Hom_{\O_X}(\S(P),\S(Q))$ defined by $u\mapsto\S(u)$ is bijective. By Lemma \ref{finitely-presented-implies-PM-is-sheaf}, $\P(Q)$ is a sheaf and hence, $\P(Q)=\S(Q)$. Therefore, \eqref{def-eq-tensor-functor}, $\Gamma\left(X,\S(Q)\right)\simeq Q$. Thus, we get the following isomorphism, and hence the result.
	\[\Hom_A(P,Q)\xrightarrow{\sim}\Hom_A\left(P,\Gamma\left(X,\S(Q)\right)\right)\xrightarrow{\sim}\Hom_{\O_X}\left(\S(P),\S(Q)\right)\]
\end{proof}

We recall that the cohomology of supersheaf is same as that of a cohomology of underlying sheaf of abelian groups, for more details we refer to \cite[Chapter II, Section 1]{Bruzzo91}. We also recall that a supersheaf $\F$ on a topological space $X$ is said to be \emph{acyclic} if all cohomology groups $H^k(X,\F)=0$ for all $k>0$ (see \cite[Chapter II, Definition 2.7]{Bruzzo91}).

Now we are in the position to state and prove the super analogue of the Serre-Swan theorem.

\begin{theorem}[Super analogue of the Serre-Swan theorem]\label{super-serre-swan}
	Let $(X,\O_X)$ be a locally ringed superspace and $A=\Gamma(X,\O_X)$. Assume that every locally free supersheaf of bounded rank is acyclic and generated by finitely many global sections. Then the functor $\S$ defines an equivalence of categories from the category $\SFgp{A}$ of finitely generated superprojective right $A$-supermodules to the category $\Slfb{X}$ of locally free supersheaves of bounded rank over $X$.
\end{theorem}
\begin{proof}
	As $\S$ is fully faithful by Proposition \ref{S is fully faithful}, it is enough to show that $\S$ is essentially surjective, i.e., for a locally free supersheaf $\F$ of bounded rank we need to find a finitely generated superprojective right $A$-supermodule $P$ such that $\S(P)\simeq\F$.
	
	As $\F$ is finitely generated by global sections, there exists $r,s\in\mathbb{N}$ such that $u\colon \O_X^{r|s}\to\F$ is a surjective morphism of $\O_X$-supermodules. For every in $x\in X$, we thus get an exact sequence $\O_{X,x}^{r|s}\xrightarrow{u_x}\F_x\to 0$. Note that $u_x$ is a homomorphism of $\O_{X,x}$-supermodules. As $\F$ is locally free, $\F_x$ is a free $\O_{X,x}$-supermodule. Hence, the following sequence of $\O_{X,x}$-supermodule is exact.
	\[\overline{\rm Hom}_{\O_{X,x}}\left(\F_x,\O_{X,x}^{r|s}\right)\to\overline{\rm Hom}_{\O_{X,x}}\left(\F_x,\F_x\right)\to 0\]
	By Proposition \ref{isomorphism-of-Hom-supersheaf} the following is an exact sequence.
	\[\left(\overline{\mathcal{H}om}_{\O_{X}}\left(\F,\O_{X}^{r|s}\right)\right)_x\to\left(\overline{\mathcal{H}om}_{\O_{X}}\left(\F,\F\right)\right)_x\to 0\]
	Therefore, $\overline{\mathcal{H}om}_{\O_{X}}\left(\F,\O_{X}^{r|s}\right)\to\overline{\mathcal{H}om}_{\O_{X}}\left(\F,\F\right)\to 0$ is exact. The kernel of this map is locally free by Proposition \ref{kernel-is-locally-free-of-locally-free-sheaves} and hence it is acyclic. By applying global section functor and considering the even component, we get that \[\Hom(\F,u)\colon\Hom_{\O_X}(\F,\O_X^{r|s})\to\Hom_{\O_X}(\F,\F)\] is surjective. Therefore, there exists $v\in\Hom_{\O_X}(\F,\O_X^{r|s})$ such that $\Hom(\F,u)(v)=1_\F$. Thus, $\Gamma(X,u)\colon\Gamma(X,\O_X^{r|s})\to\Gamma(X,\F)$ has a right inverse, viz., $\Gamma(X,v)$. Hence, $\Gamma(X,\F)$ is superprojective $A$-supermodule.
	We get the required result because $\S\left(\Gamma(X,\F)\right)\simeq\F$ by Proposition \ref{key-prop-to-main-theorem}.
\end{proof}

\section{Some examples}
In this section we give examples of \emph{affine superschemes} and \emph{smooth supermanifolds} where the super analogue of the Serre-Swan theorem \ref{super-serre-swan} holds.
\subsection{Affine superschemes}
\begin{definition}[Prime spectrum of a supercommutative ring]\label{Prime spectrum of a supercommutative ring}
	Let $A$ be a supercommutative ring such that $A=A_0\bigoplus A_1$. The \emph{prime spectrum of a supercommutative ring} $A$, denoted by $\spec A$, is the topological space of prime ideals of the commutative ring $A_0$, i.e., $\spec A=\spec A_0$. We refer to \cite[Section 5.4]{Westra} for more details.
\end{definition}

\begin{remark}\label{basis for spec}
	Recall that for an element $f\in A$, the \emph{principle open subset} of $\spec A$ is the following. \[D(f)=\left\{\p\in\spec A:f\not\in\p\right\}\]
	Note that if $f=f_0+f_1$ then $f_1$ belongs to every prime ideal of $A$, and hence $D(f)=D(f_0)$. So it is enough to consider principle open sets and localization at even elements.
\end{remark}

\begin{definition}[Supersheaf associated to a module]\label{Supersheaf associated to module}
	Let $A$ be a supercommutative ring and $M$ a right $A$-supermodule. The sheaf $\widetilde{M}$ associated to $M$ is such that $\widetilde{M}(D(f))=M_f\simeq M\otimes_A A_f$, for an even element $f\in A_0$. For more details we refer to \cite[Section 6.4]{Westra}
\end{definition}

\begin{theorem}[Acyclicity of affine supersheaves]\label{Acyclicity of affine supersheaves}
	Let $(X,\O_X)$ be an affine superscheme and let $\F$ be a quasi-coherent supersheaf of $\O_X$-supermodules on $X$. Then $H^p(X,\F)=0$ for all $p\geq 1$.
\end{theorem}

\begin{remark}
	A proof of the above theorem follows from J.-P. Serre's result \cite[Section 44]{Serre} and one of the Cartan's result on isomorphism of \v{C}ech and sheaf cohomology groups, see \cite[Chapter III. Corollary 4 to Theorem 3.8.5]{Grothendieck}.
\end{remark}

\begin{proposition}
	Let $A$ be a superring. Theorem \ref{super-serre-swan} holds for an affine superscheme $(\spec A,\widetilde{A})$.
\end{proposition}
\begin{proof}
	Let $\F$ be a locally free supersheaf of bounded rank. Hence, $\F$ is quasi-coherent and by Theorem \ref{Acyclicity of affine supersheaves}, $\F$ is acyclic. By \cite[II.5.5]{Hartshorne}, there is an equivalence of categories between quasi-coherent sheaves on $\spec A$ and the category of $A$-modules. Note that the same proof carries over for the supermodule structures because of Remark \ref{basis for spec}. So there exists a right $A$-supermodule $M$ such that $\F\simeq\widetilde{M}$. Note that $\widetilde{M}$ is generated by global sections and hence so is $\F$. Since $\F$ is locally free of finite rank and $\spec A$ is quasicompact, $\F\simeq\widetilde{M}$ is finitely generated by global sections.
\end{proof}

\subsection{Smooth supermanifolds}
We recall the definition of the \emph{fine sheaf}. We assume that paracompact spaces are Hausdorff.
\begin{definition}[Fine sheaf]\label{def-fine-sheaf}
	Let $X$ be a paracompact topological space and let $\F$ be a sheaf of abelian groups on $X$. The sheaf $\F$ is called a \emph{fine sheaf} if for any locally finite open cover $\{U_i\}_{i\in I}$ of $X$, there exists a family of endomorphisms $\{\eta_i\}_{i\in I}$ of $\F$ such that 
	\begin{enumerate}
		\item ${\rm supp}(\eta_i)\subset U_i$, where ${\rm supp}(\eta_i)=\overline{\{x\in X:(\eta_i)_x\neq 0\}}$;
		\item $\sum\eta_i=1_\F$.
	\end{enumerate}
	The family $\{\eta_i\}_{i\in I}$ is called a \emph{partition of unity} of $\F$ subordinate to a covering of $\{U_i\}_{i\in I}$.
\end{definition}
\begin{remark}
	We keep the notations of Definition \ref{def-fine-sheaf}. Suppose that $\F$ is a fine supersheaf. We claim that a partition of unity $\{\eta_i\}_{i\in I}$ of $\F$ can be chosen so that $\eta_i\colon\F\to\F$ is an endomorphism of supersheaf for each $i\in I$. Indeed, consider a partition of unity $\{\zeta_i\}_{i\in I}$ of $\F$. We write $\zeta_i=\zeta_{0i}+\zeta_{1i}$ for each $i\in I$, where $\zeta_{0i}$ is an even endomorphism of $\F$ and $\zeta_{1i}$ is odd. Note that ${\rm supp}(\zeta_{0i})\subseteq{\rm supp}(\zeta_i)\subset U_i$. Furthermore, as $\sum\zeta_i=1$, we have that $\sum\zeta_{0i}+\sum\zeta_{1i}=1$. As the identity is an even endomorphism, we get $\sum\zeta_{1i}=0$. So, we may consider $\{\zeta_{0i}\}_{i\in I}$ to be the required partition of unity.
\end{remark}

\begin{example}
	The supersheaf of continuous functions $\mathcal{C}_\Lambda$ as defined in Example \ref{supersheaf-of-continous-functions} is a fine sheaf when $X$ is a paracompact topological space. Since $X$ is paracompact space, it admits a partition of unity, say $\{f_i\}_{i\in I}$ subordinate to (locally finite) open covering $\{U_i\}_{i\in I}$. We define $\eta_i\colon\mathcal{C}_\Lambda\to\mathcal{C}_\Lambda$ as follows.
	
	For every open subset $U\subseteq X$, $(\eta_i)_U\colon\mathcal{C}_\Lambda(U)\to\mathcal{C}_\Lambda(U)$ is defined by $f\mapsto f_i|_U\cdot f$. Note that ${\rm supp}(\eta_i)={\rm supp}(f_i)$ and $\sum\eta_i=1_{\mathcal{C}_\Lambda}$.
	
	Similarly the sheaf $\mathcal{C}_\Lambda^\infty$ is a fine sheaf.
\end{example}

We have the following result from \cite[Proposition 2.2.2]{Archana}.
\begin{proposition}\label{fine-sheaf-generation-by-global-sections}
	Let $(X,\O_X)$ be a ringed superspace such that $X$ is a paracompact topological space, and let $\O_X$ be a fine supersheaf. Consider an $\O_X$-module $\F$. Let $x\in X$, $U$ be an open neighborhood of $x$, and let $s'\in\F(U)$. Then, there exists an open neighborhood $V$ of $x\in U$, and a global section $s$ of $\F$ such that $s|_V=s'|_V$. In particular, the canonical homomorphism $\rho_x\colon\Gamma(X,\F)\to\F_x$ is surjective, and hence $\F$ is generated by global sections.
\end{proposition}

\begin{corollary}
	Let $(X,\mathcal{C}_\Lambda)$ be as in Example \ref{supersheaf-of-continous-functions} with $X$ a paracompact topological space. Then $\Gamma(X,\cdot)\colon{\bf mod}\text{- }\mathcal{C}_\Lambda\to{\bf mod}\text{- }\Gamma(X,\mathcal{C}_\Lambda)$ is fully faithful functor.
\end{corollary}

\begin{theorem}\label{serre-swan-for-fine-sheaves}
	Let $(X,\O_X)$ locally ringed superspace such that $X$ is compact and $\O_X$ a fine supersheaf. Then Theorem \ref{super-serre-swan} holds.
\end{theorem}
\begin{proof}
	Let $\F$ be a locally free $\O_X$-supermodule of bounded rank. By Proposition \ref{fine-sheaf-generation-by-global-sections}, for any $x\in X$ there is a finite set $S_x=\{s_{1},\ldots,s_{n}\}\subset\Gamma(X,\F)$ such that $\{(s_i)_x\}$ generates $\F_x$. Hence, for any $x\in X$ there exists a neighborhood $U_x$ of $x$ such that for any $y\in U_x$, $\{(s_i)_y\}$ generates $\F_y$. As $X$ is compact there are finitely many $x_i$ such that $X=\bigcup U_{x_i}$ and thus, for any $x\in X$, the germs of $S=\bigcup S_{x_i}$ will generate $\F_x$. This shows that $\F$ is generated by finitely many global sections. Since, $\O_X$ is a fine supersheaf so is $\F$. By \cite[II.3.5]{Wells}, $\F$ is acyclic. Hence the result.
\end{proof}

For the definition of a \emph{smooth supermanifold} we follow \cite{Batchelor}.

\begin{definition}[Smooth supermanifold]\label{def-supermanifold}
	Let $X$ be a second countable, Hausdorff topological space and let $(X,\O_X)$ be a ringed superspace, where $\O_X$ is a supersheaf of $\R$-superalgebras. The $(X,\O_X)$ is a \emph{smooth supermanifold} of \emph{odd dimension} $n$, if the following conditions are satisfied.
	\begin{itemize}
		\item There is a morphism of supersheaves $\O_X\to C^\infty$, where the grading on $C^\infty$ is given by $\left(C^\infty(U)\right)_0=C^\infty(U)$ for any open subset $U\subseteq X$.
		\item (Local triviality conditions). There exists an open cover $\{U_i\}_{i\in I}$ of $X$ such that there are isomorphisms of supersheaves $\O_X|_{U_i}\simeq C^\infty|_{U_i}\otimes_{\O_X|_{U_i}}\bigwedge\R^n$ for every $i\in I$.
	\end{itemize}
\end{definition}

\begin{proposition}\label{compact supermanifold fine sheaf}
	If $(X,\O_X)$ be a compact smooth supermanifold then, $\O_X$ is a fine supersheaf.
\end{proposition}
\begin{proof}
	We fix isomorphisms $\phi_i\colon\O_X|_{U_i}\to C^\infty|_{U_i}\otimes_{\O_X|_{U_i}}\bigwedge\R^n$ in the local triviality conditions (see Definition \ref{def-supermanifold}). Let $\{V_j\}_{j\in J}$ be a locally finite open cover of $X$. Since $X$ is assumed to be compact, there exists a smooth partition unity subordinate to $V_j$, viz., $f_j\colon X\to\R$ for every $j\in J$.  We now find a partition of unity $\eta_j\colon\O_X\to\O_X$. In this regard, we define for each $x\in X$, a family of maps $\eta_j^x\colon\O_{X,x}\to\O_{X,x}$.
	
	Let $x\in X$. There exists $i\in I$ such that $x\in U_i$. We define $\eta_{j}^x\colon\O_{X,x}\to\O_{X,x}$ by $\eta_{j}^x(s_x)=\phi^{-1}_{i,x}\left(\phi_{i,x}(s_x)\cdot f_{j,x}\right)$. We now consider an open subset $U\subseteq X$ and a section $s\in\O_X(U)$. We show that for all $y\in U$ there exists open neighborhood $V_y\subseteq U$ of $y$ and $t\in\O_X(V_y)$ such that $t_y=\eta_j^x(s_y)$. Indeed, we can take $V_y=U_k\cap U$, where $y\in U_k$, and $t=(\phi_k)_{(U_k\cap U)}^{-1}\left((\phi_k)_{(U_k\cap U)}\left(s|_{U_k\cap U}\right)\cdot f_k\right)$. By Proposition \ref{local-morphism-to-global-morphism}, there exists $\eta_j\colon\O_X\to\O_X$. It follows that ${\rm supp}(\eta_j)\subseteq V_j$ and $\sum\eta_j=1$, i.e., the family $\{\eta_j\}$ is a partition of unity. This shows that $\O_X$ is a fine supersheaf.
\end{proof}

We obtain the following corollary to Theorem \ref{serre-swan-for-fine-sheaves} by using the above Proposition \ref{compact supermanifold fine sheaf}.

\begin{corollary}
	Let $(X,\O_X)$ be a smooth supermanifold and assume that $X$ is compact. Then, Theorem \ref{super-serre-swan} holds.
\end{corollary}

\section*{Acknowledgment}
We thank N. Raghavendra for useful suggestions. We are also grateful to the referee for their careful reading and constructive comments which helped improve this paper.

\bibliographystyle{plain}
\bibliography{super}

@book{MacLane,
  author = {Mac Lane, Saunders},
  title = {Categories for the working mathematician.},
  edition = {2nd ed},
  fseries = {Graduate Texts in Mathematics},
  series = {Grad. Texts Math.},
  issn = {0072-5285},
  volume = {5},
  isbn = {0-387-98403-8},
  year = {1998},
  publisher = {New York, NY: Springer},
  language = {English},
  doi = {book/10.1007/978-1-4757-4721-8},
  zbMATH = {1216133},
  Zbl = {0906.18001}
 }

@book{Bourbaki-commutative,
  author = {Bourbaki, Nicolas},
  title = {Elements of mathematics. {Commutative} algebra. {Chapters} 1--7. {Transl}. from the {French}.},
  edition = {Softcover edition of the 2nd printing 1989},
  isbn = {3-540-64239-0},
  year = {1998},
  publisher = {Berlin: Springer},
  language = {English},
  zbMATH = {1194716},
  Zbl = {0902.13001}
 }

@article{Morye-2013,
  author = {Morye, Archana S.},
  title = {Note on the {Serre}-{Swan} theorem},
  fjournal = {Mathematische Nachrichten},
  journal = {Math. Nachr.},
  issn = {0025-584X},
  volume = {286},
  number = {2-3},
  pages = {272--278},
  year = {2013},
  language = {English},
  doi = {10.1002/mana.200810263},
  zbMATH = {6145346},
  Zbl = {1274.14020}
 }

@article {Batchelor,
  AUTHOR = {Batchelor, Marjorie},
  TITLE = {The structure of supermanifolds},
  JOURNAL = {Trans. Amer. Math. Soc.},
  FJOURNAL = {Transactions of the American Mathematical Society},
  VOLUME = {253},
  YEAR = {1979},
  PAGES = {329--338},
  ISSN = {0002-9947,1088-6850},
  MRCLASS = {58A05 (83E99)},
  MRNUMBER = {536951},
  MRREVIEWER = {A.\ Verona},
  DOI = {10.2307/1998201},
  URL = {https://doi.org/10.2307/1998201},
 }

@article {Grothendieck,
  AUTHOR = {Grothendieck, Alexander},
  TITLE = {Sur quelques points d'alg\`ebre homologique},
  JOURNAL = {Tohoku Math. J. (2)},
  FJOURNAL = {The Tohoku Mathematical Journal. Second Series},
  VOLUME = {9},
  YEAR = {1957},
  PAGES = {119--221},
  ISSN = {0040-8735,2186-585X},
  MRCLASS = {18.00},
  MRNUMBER = {102537},
  MRREVIEWER = {D.\ Buchsbaum},
  DOI = {10.2748/tmj/1178244839},
  URL = {https://doi.org/10.2748/tmj/1178244839},
 }

@book {Hartshorne,
  AUTHOR = {Hartshorne, Robin},
  TITLE = {Algebraic geometry},
  SERIES = {Graduate Texts in Mathematics},
  VOLUME = {No. 52},
  PUBLISHER = {Springer-Verlag, New York-Heidelberg},
  YEAR = {1977},
  PAGES = {xvi+496},
  ISBN = {0-387-90244-9},
  MRCLASS = {14-01},
  MRNUMBER = {463157},
  MRREVIEWER = {Robert\ Speiser},
 }

@article {MSV,
  AUTHOR = {Morye, Archana S. and Sarma Phukon, Aditya and Devichandrika,
   V.},
  TITLE = {Notes on super projective modules},
  JOURNAL = {Indian J. Pure Appl. Math.},
  FJOURNAL = {Indian Journal of Pure and Applied Mathematics},
  VOLUME = {54},
  YEAR = {2023},
  NUMBER = {4},
  PAGES = {1226--1238},
  ISSN = {0019-5588,0975-7465},
  MRCLASS = {16W55 (16D40 16W50)},
  MRNUMBER = {4662837},
  MRREVIEWER = {M.\ A.\ Raza},
  DOI = {10.1007/s13226-022-00336-4},
  URL = {https://doi.org/10.1007/s13226-022-00336-4},
 }

@article {Serre,
  AUTHOR = {Serre, Jean-Pierre},
  TITLE = {Faisceaux alg\'ebriques coh\'erents},
  JOURNAL = {Ann. of Math. (2)},
  FJOURNAL = {Annals of Mathematics. Second Series},
  VOLUME = {61},
  YEAR = {1955},
  PAGES = {197--278},
  ISSN = {0003-486X},
  MRCLASS = {14.0X},
  MRNUMBER = {68874},
  MRREVIEWER = {C.\ Chevalley},
  DOI = {10.2307/1969915},
  URL = {https://doi.org/10.2307/1969915},
 }

@book {Wells,
  AUTHOR = {Wells, Jr., R. O.},
  TITLE = {Differential analysis on complex manifolds},
  SERIES = {Graduate Texts in Mathematics},
  VOLUME = {65},
  EDITION = {Second},
  PUBLISHER = {Springer-Verlag, New York-Berlin},
  YEAR = {1980},
  PAGES = {x+260},
  ISBN = {0-387-90419-0},
  MRCLASS = {58-01 (32-01)},
  MRNUMBER = {608414},
  MRREVIEWER = {Daniel\ M.\ Burns, Jr.},
 }

@book {Bruzzo91,
  AUTHOR = {Bartocci, Claudio and Bruzzo, Ugo and Hern\'andez Ruip\'erez,
   Daniel},
  TITLE = {The geometry of supermanifolds},
  SERIES = {Mathematics and its Applications},
  VOLUME = {71},
  PUBLISHER = {Kluwer Academic Publishers Group, Dordrecht},
  YEAR = {1991},
  PAGES = {xx+242},
  ISBN = {0-7923-1440-9},
  MRCLASS = {58A50 (32C11 58C50)},
  MRNUMBER = {1175751},
  MRREVIEWER = {Jaime\ Mu\~noz Masqu\'e},
  DOI = {10.1007/978-94-011-3504-7},
  URL = {https://doi.org/10.1007/978-94-011-3504-7},
 }

@article {Bruzzo,
  AUTHOR = {Bruzzo, Ugo and Hern\'andez Ruip\'erez, Daniel and Polishchuk,
   Alexander},
  TITLE = {Notes on fundamental algebraic supergeometry. {H}ilbert and
   {P}icard superschemes},
  JOURNAL = {Adv. Math.},
  FJOURNAL = {Advances in Mathematics},
  VOLUME = {415},
  YEAR = {2023},
  PAGES = {Paper No. 108890, 115},
  ISSN = {0001-8708,1090-2082},
  MRCLASS = {14M30 (14D22 14H10 14K10 83E30)},
  MRNUMBER = {4544562},
  MRREVIEWER = {Frans\ Cantrijn},
  DOI = {10.1016/j.aim.2023.108890},
  URL = {https://doi.org/10.1016/j.aim.2023.108890},
 }

@book {Manin,
  AUTHOR = {Manin, Yuri I.},
  TITLE = {Gauge field theory and complex geometry},
  SERIES = {Grundlehren der mathematischen Wissenschaften [Fundamental
   Principles of Mathematical Sciences]},
  VOLUME = {289},
  NOTE = {Translated from the Russian by N. Koblitz and J. R. King},
  PUBLISHER = {Springer-Verlag, Berlin},
  YEAR = {1988},
  PAGES = {x+297},
  ISBN = {3-540-18275-6},
  MRCLASS = {32-02 (14-02 32L25 58A50 81-02)},
  MRNUMBER = {954833},
 }

@book {Archana,
  AUTHOR = {Morye, Archana S.},
  TITLE = { On the Serre-Swan theorem, and on vector bundles over Real abelian varieties},
  NOTE = {Thesis (Ph.D.)--Harish-Chandra Research Institute, Allahabad},
  PUBLISHER = {},
  YEAR = {2011},
  PAGES = {127},
  ISBN = {},
  MRCLASS = {},
  MRNUMBER = {},
  URL =
  {https://www.hri.res.in/~libweb/theses/softcopy/thes_archana_morye.pdf},
 }

@book {Westra,
  AUTHOR = {Westra, Dennis B.},
  TITLE = {Superrings and Supergroups},
  NOTE = {Thesis (Ph.D.)--Universit\"at of Wien},
  PUBLISHER = {},
  YEAR = {2009},
  PAGES = {190},
  ISBN = {},
  MRCLASS = {},
  MRNUMBER = {},
  URL =
  {https://www.mat.univie.ac.at/~michor/westra_diss.pdf},
 }

@article {Swan,
  AUTHOR = {Swan, Richard G.},
  TITLE = {Vector bundles and projective modules},
  JOURNAL = {Trans. Amer. Math. Soc.},
  FJOURNAL = {Transactions of the American Mathematical Society},
  VOLUME = {105},
  YEAR = {1962},
  PAGES = {264--277},
  ISSN = {0002-9947,1088-6850},
  MRCLASS = {57.30 (13.40)},
  MRNUMBER = {143225},
  MRREVIEWER = {K.\ Kondo},
  DOI = {10.2307/1993627},
  URL = {https://doi.org/10.2307/1993627},
 }

	\end{document}